\documentclass[a4paper,12pt]{article}
\usepackage{graphicx} \usepackage{float}
\usepackage{amsmath} \usepackage{amssymb} \usepackage{amsfonts}
\usepackage{xcolor} \usepackage{stmaryrd}
\usepackage{array}
\usepackage{booktabs}
\usepackage{multirow}

\usepackage[pdftex]{hyperref}
\usepackage[round, comma]{natbib}

\usepackage[utf8]{inputenc}
\usepackage[T1]{fontenc}
\usepackage{lmodern}
\usepackage{algorithm2e}
\numberwithin{equation}{section}

\newtheorem{rem}{Remark}[section]
\newtheorem{theo}{Theorem}[section]

\newtheorem{pro}{Proposition}[section]

\newcommand{\N}{\mathbb{N}}
\newcommand{\R}{\mathbb{R}}

\newcommand{\norme}[2]{{\left\|#2\right\|_{#1}}}
\newcommand{\pnorme}[2]{{\|#2\|_{#1}}}
\newcommand{\bnorme}[2]{{\big\|#2\big\|_{#1}}}
\newcommand{\op}[1]{\mathop{\mathrm{#1}}}
\newcommand{\gtau}{\mathcal O}

\newcommand{\bZ}{\mathbf{Z}}

\newcommand{\e}{\mathbb{E}}
\newcommand{\var}{\mathbb{V}}
\newcommand{\p}{\mathbb{P}}
\newcommand{\1}{\mathbb{I}}
\newcommand{\given}{|}

\newcommand{\te}{\widetilde{\e}}
\newcommand{\tte}{\widetilde{\widetilde\e}}

\renewcommand{\hat}{\widehat}

\newcommand{\interval}[1]{{\llbracket1,#1\rrbracket}}
\newcommand{\ffi}{{f_i}}
\newcommand{\bfi}{{\bar f_{i}}}
\newcommand{\hfi}{{\hat f_{i}}}
\newcommand{\Sin}{{S_{i,n}}}
\newcommand{\hNi}{{\hat N_i}}

\newcommand{\Nj}{{N_i}}
\newcommand{\gi}{{g_{i,n}}}
\newcommand{\ai}{{\alpha_i}}
\newcommand{\bai}{{\bar\alpha_i}}
\newcommand{\hai}{{\hat\alpha_i}}
\newcommand{\cluster}[1]{{\mathcal{X}_{#1}}}

\newcounter{mnotecount}[section]

\graphicspath{{FIGURES/}}
\hypersetup{%
  breaklinks=true,
  colorlinks=true,
}

\begin{document}
\begin{center}

{\sc On clustering procedures and nonparametric mixture estimation\\
\vspace{0.7cm}}

Stéphane Auray $^{\mbox{\footnotesize a}}$, Nicolas Klutchnikoff $^{\mbox{\footnotesize b}}$ and
Laurent Rouvière $^{\mbox{\footnotesize c,}}\footnote{Corresponding
  author.}^{,}$

\vspace{0.5cm}

$^{\mbox{\footnotesize a}}$ CREST-Ensai \\
 EQUIPPE (EA4018) -- ULCO\\ and CIRPEE, Canada\\
\smallskip
\textsf{stephane.auray@ensai.fr}\\
\bigskip
$^{\mbox{\footnotesize b}}$ CREST-Ensai, Institut de recherche mathématique avancée\\
(UMR 7501, CNRS and Universit\'e de Strasbourg)\\
\smallskip
\textsf{nicolas.klutchnikoff@ensai.fr}
\bigskip

$^{\mbox{\footnotesize c}}$ CREST-Ensai, IRMAR (UMR 6625), UEB
\smallskip

\textsf{laurent.rouviere@ensai.fr}\\
\vspace{0.5cm}

$^{\mbox{\footnotesize a,b,c}}$ Campus de Ker-Lann, Rue Blaise Pascal, BP 37203\\
35172 Bruz cedex, France

\end{center}

\begin{abstract}
This paper deals with  nonparametric estimation of conditional densities in mixture models in the case when additional covariates are available. The proposed approach consists of performing a preliminary clustering algorithm on the additional covariates to guess the mixture component of each observation. Conditional densities of the mixture model are then estimated using kernel density estimates applied separately to each cluster. We investigate the expected $L_1$-error of the resulting estimates and derive optimal rates of convergence over classical nonparametric density classes provided the clustering method is accurate. Performances of clustering algorithms are measured by the \emph{maximal misclassification error}. We obtain upper bounds of this quantity for a single linkage hierarchical clustering algorithm. Lastly, applications of the proposed method to mixture models involving electricity distribution data and simulated data are presented.
\medskip

\footnotesize{\noindent
}
\end{abstract}

\section{Introduction}
Finite mixture models are widely used to account for population heterogeneities. In many fields such as biology, econometrics and social sciences, experiments are based on the analysis of a variable characterized by a different behavior depending on the group of individuals. A natural way to model heterogeneity for a real random variable $Y$ is to use a mixture model. In this case, the density $f$ of $Y$ can be written as
\begin{equation}
\label{eq:defmixtY}
f(t)=\sum_{i=1}^M\ai\ffi(t),\quad t\in\R.
\end{equation}
Here $M$ is the number of subpopulations, $\ai$ and $\ffi$ are respectively the mixture proportion and the probability density function of the $i$\textsuperscript{th} subpopulation. We refer the reader to \cite{evehan81,macbas88,macpeel00} for a broader picture of mixture density models as well as for practical applications. \bigskip

When dealing with mixture density models such as~\eqref{eq:defmixtY}, some issues arise. In some cases, the number of components $M$ is unknown and needs to be estimated. To this end, some algorithms have been developed to provide consistent estimates of this parameter. For instance, when $M$ corresponds to the number of modes of $f$, \cite{cufefr00} and~\cite{MR2320821} propose an estimator based on the level sets of $f$. Model identifiability is an additional issue that has received some attention in the literature. Actually, model~\eqref{eq:defmixtY} is identifiable only by imposing restrictions on the vector $(\alpha_1,\hdots,\alpha_M,f_1,\hdots,f_M)$. In order to provide the minimal assumptions such that ~\eqref{eq:defmixtY} becomes identifiable, \cite{celgov95,bomova06} (see also the references therein) assume that the density functions $\ffi$'s belong to some parametric or semi-parametric density families. However, in a nonparametric setting, it turns out that identifiability conditions are more difficult to provide. \cite{halzho03} define mild regularity conditions to achieve identifiability in a multivariate nonparametric setting while~\cite{kita04} considers the case where appropriate covariates are available. \bigskip

When the model~\eqref{eq:defmixtY} is identifiable, the statistical problem consists of estimating mixture proportions $\ai$ and density functions $\ffi$. In the parametric case, some algorithms have been proposed such as maximum likelihood techniques (\cite{lind83a,lind83b,redwal84}) as well as Bayesian approaches (\cite{dierob94,bicego00}). When the $\ffi$'s belong to nonparametric families, it is often assumed that training data are observed, {\it i.e.}, the component of the mixture from which $Y$ is distributed is available. In that case, the model is identifiable and some algorithms allow to estimate both the $\ai$'s and the $\ffi$'s (see~\cite{titt83,haltit84,haltit85,cerr92}). However, as pointed out by~\cite{halzho03}, inference in mixture nonparametric density models becomes more difficult without training data. These authors introduce consistent nonparametric estimators of the conditional distributions in a multivariate setting. We also refer to \cite{bomova06} who provide efficient estimators under the assumption that the unknown mixed distribution is symmetric. These estimates are extended by~\cite{bechhu09,bechhu11} for multivariate mixture models.

\bigskip

The framework we consider takes place between the two above situations. More precisely, training data are not observed but we assume to have at hand some covariates that may provide information on the components of the mixture from which $Y$ is distributed. Our approach consists of performing a preliminary clustering algorithm on these covariates to guess the mixture component of each observation. Density functions $\ffi$ are then estimated using a nonparametric density estimate based on the predictions of the clustering method. \bigskip

Many authors have already proposed to carry out a preliminary clustering step to improve density estimates in mixture models. \cite{ruruka06} conduct a comprehensive simulation study to conclude that
a preliminary clustering using the EM algorithm allows to some extent to improve performances of some density estimates (see also \cite{jeolan94}). However, to our knowledge, no work has been devoted so far to measure the effects of the clustering algorithm on the resulting estimates of the distribution functions $\ffi$. This paper proposes to fill this gap, studying the $L_1$-error of these estimates. To do so, we measure the performance of clustering methods by the maximal misclassification error \eqref{eq:def_phin}. This criterion allows us to derive optimal rates of convergence over classical nonparametric density classes, provided the clustering method used in the first step performs well with respect to this notion.
\bigskip

The paper is organized as follows. In Section \ref{sec:model}, we present the two-step estimator and give the main results. Examples of clustering algorithms are worked out in Section \ref{sec:cluster}. In particular, the maximal misclassification error of a hierarchical clustering algorithm is studied under mild assumptions on the model. Applications on simulated and real data are presented in Sections \ref{sec:simu} and \ref{sec:erdf}. A short conclusion including a discussion of the implications of the work is given in Section \ref{sec:conclusion} and proofs are gathered in Section \ref{sec:proofs}.

\section{A two-step nonparametric estimator}
\label{sec:model}

\subsection{The statistical problem}
\label{sec:presentation}
Our focus is on the estimation of conditional densities in a univariate mixture density model. Formally we let $(Y,I)$ be a random vector taking values in $\R\times \interval{M}$ where $M\geq
2$ is a known integer. We assume that the distribution of $Y$ is
characterized by a density $f$ defined, for all $t\in\R$, by
\begin{equation*}
f(t)=\sum_{i=1}^M \ai \ffi(t),
\end{equation*}
where, for all $i\in\interval{M}$, $\ai=\p(I=i)$ are the prior probabilities (or the weights of the mixture) and $\ffi$ are the densities of the conditional distributions $\mathcal L(Y|I=i)$ (or
the components of the mixture).
\bigskip

If we have at hand $n$ observations $(Y_1,I_1),\hdots,(Y_n,I_n)$ drawn from the distribution of $(Y,I)$, one can easily find efficient estimates for both the $\ai$'s and the $\ffi$'s. For example, if
we denote $\Nj = \#\left\{k\in\interval{n} : I_k=i\right\}$, then we can estimate $\ai$ using the empirical proportion $\bai=\Nj/n$ and $\ffi$ by the kernel density estimate $\bfi$ defined for all $t\in\R$ by
\begin{equation}
\label{eq:defbarf}
\bfi(t) = \frac1{\Nj} \sum_{k=1}^n K_h(t,Y_k) \1_i(I_k)
\end{equation}
if $N_i>0$. For the definiteness of $\bfi$ we conventionally set $\bfi(t)=0$ if $N_i=0$.
Here $K$ is a kernel which belongs to $L_1(\R,\R)$ and such that $\int K=1$, $h>0$ is a bandwidth and
\begin{equation}
\label{eq:def_Kth}
K_{h}(t,y)=\frac{1}{h}K\left(\frac{t-y}{h}\right)
\end{equation}
is the classical convolution kernel located at point $t$ (see \cite{ros56} and~\cite{par62} for instance). Estimate \eqref{eq:defbarf} is just the usual kernel density estimate defined from observations in the $i^{\textrm{th}}$ subpopulation. It follows that, under classical assumptions regarding the smoothing parameter $h$ and the kernel $K$, $\bfi$ has similar properties as those of the well-known kernel density estimate. In particular, the expected $L_1$-error
$$\e\|\bfi-\ffi\|_1=\e\int_\R|\bfi(t)-\ffi(t)|dt$$
achieves optimal rates when $\ffi$ belongs to regular  density classes such as Hölder or Lipschitz classes (see~\cite{devgyo85}).
\bigskip

The problem is more complicated when the random variable $I$ is not observed. In this situation,  $\bai$ and $\bfi$ are not computable and one has to find another way to define efficient estimates for both $\ai$ and $\ffi$. In this work, we assume that one can obtain information on $I$ through another covariate $X$ which takes values in $\R^d$ where $d\geq 1$. This random variable is observed and its conditional distribution $\mathcal{L}(X\given I=i)$ is characterized by a density $g_i=\gi:\R^d\to\R$ which could depend on $n$. In this framework, the statistical problem is to estimate both the components and the weights of the mixture model \eqref{eq:defmixtY} using the $n$-sample $(Y_1, X_1),\ldots,(Y_n, X_n)$ extracted from $(Y_1, X_1, I_1), \ldots,(Y_n, X_n, I_n)$ randomly drawn from the distribution of $(Y,X,I)$.

\subsection{Discussion on the model}
Estimating components of a mixture model is a classical statistical problem. The new feature proposed here is to include covariates in the model which can potentially improve traditional algorithms. These covariates are represented by a random vector $X$ which provides information on the unobserved group $I$. This model includes many practical situations. Three examples are provided in this section.

\paragraph{The classical mixture problem without covariates.}
A traditional problem in mixture models is the  estimation of the components $f_i,i\in\interval{M}$ in \eqref{eq:defmixtY} from (only) an i.i.d sample $Y_1,\hdots,Y_n$ drawn from $f$: no covariates are available. In this context, many parametric methods such as the EM algorithm (and its derivatives) as well as nonparametric procedures (under suitable identifiability constraints) can be used and are widely studied. Even if this model is formally a particular case of ours (we just have to take $X=Y$), the approach presented in this paper is not designed to be competitive in this situation with dedicated parametric or nonparametric methods. Indeed, our model focus on practical situations where covariates can be used to obtain useful information about the hidden variable $I$. Below, we offer two realistic situations where such covariates are naturally available.

\paragraph{Medical example.}
Many diseases evolve over time and exhibit different stages of development which can be represented by a variable $I$ that takes a finite number of values. In many situations, the problem is not to study the stage $I$ but some variables that can potentially have different behavior according to $I$. For instance, the survival time $Y$ and its conditional distributions with respect to $I$ are typically of interest in many situations. In practice, the stage $I$ is generally not observed. It is assessed by the medical team from several items such as physiological data, medical examinations, interviews with the patient (and so on) that can be represented by covariates $X$ in our model.

\paragraph{Electricity distribution.}
A distribution network may locally experience minor problems, due for example to bad weather, that may affect some customers during a fixed period of time in a given geographical area. To better understand the origin and/or consequences of the dysfunctions, and thus better forecast network operations, electricity distributors are interested in the distribution behavior of several quantities $Y$ for two different groups of customers: those affected by the malfunction and the others. Variables $Y$ may for instance represent averages or variations of consumption after the disruption period. In this situation the group is represented by a variable $I$: $I=1$ for the users affected by the disruption and 2 for the others. This binary variable $I$ is not directly observed but it can be guessed from individuals curves of consumptions during the disruption period. In our framework, discrete versions of these curves correspond to the covariate $X$. This example is explained in-depth and analyzed in Section~\ref{sec:erdf} using real data from ERDF, the main French distributor of electricity.

\subsection{A kernel density estimate based on a clustering approach}
To estimate densities $f_i$ of the conditional distributions $\mathcal L(Y|I=i),i\in\interval{M}$, we propose a two-step algorithm that can be summarized as follows.
\begin{enumerate}
\item Apply a clustering algorithm on the sample $X_1,\hdots,X_n$ to predict the label $I_k$ of each observation $X_k$;
\item Estimate conditional densities $\ffi$ by kernel density estimates \eqref{eq:defbarf} where unobserved labels are substituted by predicted labels.
\end{enumerate}
Formally, we first perform a given clustering algorithm to split the sample $X_1,\hdots,X_n$ into $M+1$ clusters $\cluster{0},\cluster{1},\hdots,\cluster{M}$ such that $\cluster{i}\neq \emptyset$ for all $i\in\interval{M}$. Clusters $\cluster{0},\cluster{1},\hdots,\cluster{M}$ satisfy
$$\bigcup_{i=0}^M\cluster{i}=\{X_1,\hdots,X_n\}\quad\textrm{and}\quad
\forall i\neq j,\ \cluster{i}\cap\cluster{j}=\emptyset.$$
We do not specify the clustering method here, some examples are discussed in Sections~\ref{sec:cluster} and~\ref{sec:simu}. Observe that we define $M+1$ clusters instead of $M$. The cluster $\cluster{0}$ (which could be empty) contains the observations for which the clustering procedure is not able to predict the label. For
example, if the clustering procedure reveals some outliers, they are collected in $\cluster{0}$ and we do not use these outliers to estimate the $\ffi$'s.
\bigskip

Once the clustering step is performed, we define the predicted labels $\hat I_k$ as
$$
\hat I_k=i\quad\textrm{if}\quad X_k\in\cluster{i},\qquad k\in\interval{n},\quad i\in\interval{M}.
$$
Observation $X_k$ is not correctly assigned to its group with probability $\p(\hat I_k\neq I_k)$. We measure the performance of the clustering algorithm by the maximal probability to not correctly attribute an observation:
\begin{equation}
\label{eq:def_phin}
\varphi_n=\max_{1\leq k\leq n}\p(\hat I_k\neq I_k).
\end{equation}
We call this error term the \emph{maximal misclassification error}. It will be studied for two clustering algorithms in Section~\ref{sec:cluster}.
\bigskip

To define our estimates, we just replace in \eqref{eq:defbarf} the true labels $I_k$ by the predicted labels $\hat I_k$. Formally, prior probabilities $\ai$ are estimated by
\begin{equation*}
  \hai = \frac{\hNi}n
\quad\text{where}\quad
  \hNi=\#\{k\in\interval{n}:\hat I_k=i\},
\end{equation*}
while for the conditional densities $\ffi$, we consider the kernel density estimator with kernel $K:\R\to\R$ and bandwidth $h>0$
\begin{equation}
\label{eq:def_kern_est}
\hfi(t) = \frac{1}{\hNi}\sum_{k: X_k\in \cluster{i}}K_{h}(t,Y_k)=\frac{1}{\hNi}\sum_{k=1}^nK_{h}(t,Y_k)\1_{\{i\}}(\hat I_k),
\end{equation}
where $K_{h}$ is defined in~\eqref{eq:def_Kth}. Observe that since for all $i\in\interval{M}$ the clusters $\cluster{i}$ are nonempty, the estimates $\hfi$ are well defined.
\medskip

Kernel estimates $\hfi$ are defined from observations in cluster $\cluster{i}$. The underlying assumption is that, for all $i\in\interval{M}$, each cluster $\cluster{i}$ collects almost all of the observations $X_k$ such that $Y_k$ is randomly drawn from $\ffi$. Under this assumption, $\varphi_n$ is expected to be small and $\hfi$ to be closed to the oracle estimates $\bfi$ defined by equation~\eqref{eq:defbarf}. This closeness is measured in the following theorem which makes the connection between the expected $L_1$-errors of $\bfi$ and $\ffi$.

\begin{theo}
\label{theo:errL1gen}
There exist positive constants $A_1-A_3$ such that, for all $n\geq 1$ and $i\in\interval{M}$
\begin{equation}
\label{eq:borne_er_L1}
\e\bnorme1{\hfi-\ffi} \leq  \e\norme1{\bfi-\ffi}+A_1\varphi_n+A_2\exp(-n)
\end{equation}
and
\begin{equation}
\label{eq:bornealphai}
\e|\hai-\ai|\leq \varphi_n+\frac{A_3}{\sqrt{n}}.
\end{equation}
\end{theo}

Constants $A_1-A_3$ are specified in the proof of the theorem. We emphasize that inequalities~\eqref{eq:borne_er_L1} and \eqref{eq:bornealphai} are non-asymptotic, that is, the bounds are valid for all $n$. If we intend to prove any consistency results regarding $\hfi$ and $\hai$, inequality~\eqref{eq:borne_er_L1} says that the maximal misclassification error $\varphi_n$ should tend to zero. Moreover, if  $\varphi_n$ tends to zero much faster than the $L_1$-error of $\bfi$, then the asymptotic performance is guaranteed to be equivalent to the one of the oracle estimate $\bfi$. The $L_1$-error of $\bfi$, with properly chosen bandwidth $h$ and kernel $K$, is known to go to zero, under standard smoothness assumptions, at rate $n^{-\frac{s}{2s+1}}$ where $s>0$ is typically an index  representing the regularity of $\ffi$. For example, when we consider Lipschitz or Hölder classes of functions with compact supports, $s$ corresponds to the number of absolutely continuous derivatives of the functions $\ffi$. In this context, if $\varphi_n=\mathcal{O}(n^{-\frac{s}{2s+1}})$, then
$$\e\bnorme1{\hfi-\ffi}=\mathcal{O}(n^{-\frac{s}{2s+1}}).$$

\begin{rem}
\label{rem:permutation}
Note that even if clusters $\cluster{1},\hdots,\cluster{M}$ are arbitrarily indexed, inequalities \eqref{eq:borne_er_L1} and \eqref{eq:bornealphai} are true whatever the choice of the indexes. However, when indexes are not chosen according to the true labels, $\varphi_n$ could be large even if the clustering procedure performs well. In this situation there exists a permutation of the indexes such that, after this permutation, the maximal misclassification error is small. More precisely it can be readily seen, using Theorem~\ref{theo:errL1gen}, that
\begin{equation}
\min_{\pi\in\Pi_M}\e\bnorme1{\hat f_{\pi(i)}-\ffi} \leq  \e\norme1{\bfi-\ffi}+A_1\min_{\pi\in\Pi_M}\varphi_n(\pi)+A_2\exp(-n)
\end{equation}
where $\Pi_M$ denotes the set of all permutations of $\interval{M}$ and $\varphi_n(\pi)$ is the maximal misclassification error of the clustering method after the permutation of the indexes:
\begin{equation}
\label{eq:phi_perm}
\varphi_n(\pi)=\max_{k=1,\hdots,n}\p(\pi(\hat I_k)\neq I_k),\quad \pi\in\Pi_M.
\end{equation}
\end{rem}

\begin{rem}
\label{rem:bw-choice}
As usual, the choice of the bandwidth $h$ reveals crucial for the performance of the kernel density estimates. However, this paper does not provide any theory to select this parameter. If automatic or adaptive  procedures are needed, they can be obtained by adjusting traditional automatic selection procedures for classical nonparametric estimators (see for example~\cite{MR1743393} or~\cite{devlug01}).
\end{rem}

\section{Clustering procedures}
\label{sec:cluster}
The proposed procedure requires a preliminary clustering algorithm performed on the sample $X_1,\hdots,X_n$. Even if any clustering algorithm could be applied in practice, it should be chosen according to the conditional distributions $\mathcal L(X|I=i),i\in\interval{M}$. More precisely, each cluster should match up with  observations drawn from one of those conditional distributions. From a theoretical point of view, for a given clustering procedure, the problem is to find upper bounds for the maximal misclassification error $\varphi_n$ to apply Theorem~\ref{theo:errL1gen}. In a parametric setting, {\it i.e.}, when conditional distributions are identified by unknown parameters, clustering algorithms are often based on efficient estimators of these unknown parameters. We provide an example in Section \ref{sec:toy_ex}. Without parametric assumptions on the distribution, the problem is more complicated. Contrary to data analysis methods such as regression or classification, there are many ways to define clustering. One of the most popular approach consists of defining clusters as  the connected components of the level sets of the density (see \cite{har75}). This amounts to saying that clusters represent high density regions of the data separated by low density regions. In this context, many authors have studied theoretical performances of clustering algorithms based on neighborhood graphs such as hierarchical or spectral clustering algorithms. In Section \ref{sec:supportdisjoint}, we extend results of \cite{mahevo09} and \cite{arr11} to our framework for a hierarchical clustering algorithm based on pairwise distances. This procedure is challenged with other clustering methods in the simulation part.

\subsection{A parametric example}
\label{sec:toy_ex}
\begin{figure}[H]
\centering
\begin{picture}(0,0)%
\includegraphics{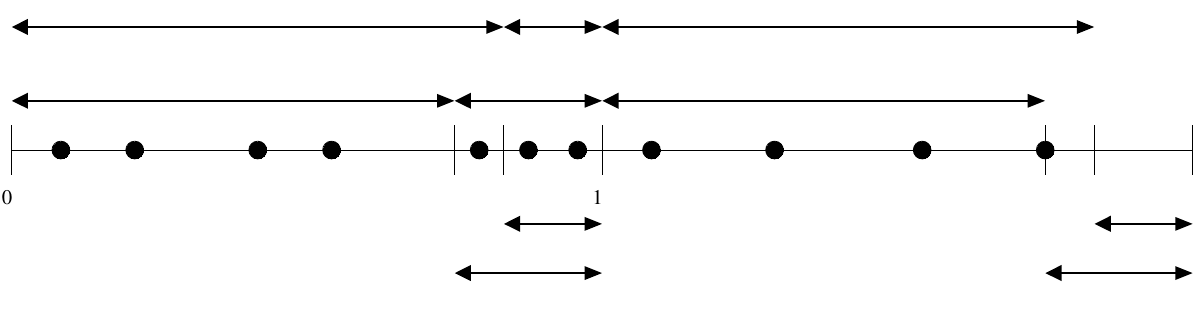}%
\end{picture}%
\setlength{\unitlength}{2072sp}%
\begingroup\makeatletter\ifx\SetFigFont\undefined%
\gdef\SetFigFont#1#2#3#4#5{%
  \reset@font\fontsize{#1}{#2pt}%
  \fontfamily{#3}\fontseries{#4}\fontshape{#5}%
  \selectfont}%
\fi\endgroup%
\begin{picture}(10927,2805)(796,-3991)
\put(2611,-1321){\makebox(0,0)[lb]{\smash{{\SetFigFont{6}{7.2}{\rmdefault}{\mddefault}{\updefault}{\color[rgb]{0,0,0}$I_k=1$}%
}}}}
\put(5311,-1321){\makebox(0,0)[lb]{\smash{{\SetFigFont{6}{7.2}{\rmdefault}{\mddefault}{\updefault}{\color[rgb]{0,0,0}$I_k=1$ or 2}%
}}}}
\put(8101,-1321){\makebox(0,0)[lb]{\smash{{\SetFigFont{6}{7.2}{\rmdefault}{\mddefault}{\updefault}{\color[rgb]{0,0,0}$I_k=2$}%
}}}}
\put(8101,-1951){\makebox(0,0)[lb]{\smash{{\SetFigFont{6}{7.2}{\rmdefault}{\mddefault}{\updefault}{\color[rgb]{0,0,0}$\hat I_k=2$}%
}}}}
\put(5221,-1951){\makebox(0,0)[lb]{\smash{{\SetFigFont{6}{7.2}{\rmdefault}{\mddefault}{\updefault}{\color[rgb]{0,0,0}$\hat I_k=0$}%
}}}}
\put(2521,-1951){\makebox(0,0)[lb]{\smash{{\SetFigFont{6}{7.2}{\rmdefault}{\mddefault}{\updefault}{\color[rgb]{0,0,0}$\hat I_k=1$}%
}}}}
\put(11611,-3031){\makebox(0,0)[lb]{\smash{{\SetFigFont{6}{7.2}{\rmdefault}{\mddefault}{\updefault}{\color[rgb]{0,0,0}2}%
}}}}
\put(5401,-3931){\makebox(0,0)[lb]{\smash{{\SetFigFont{6}{7.2}{\rmdefault}{\mddefault}{\updefault}{\color[rgb]{0,0,0}$\hat\lambda_n$}%
}}}}
\put(5761,-3391){\makebox(0,0)[lb]{\smash{{\SetFigFont{6}{7.2}{\rmdefault}{\mddefault}{\updefault}{\color[rgb]{0,0,0}$\lambda_n$}%
}}}}
\put(11161,-3391){\makebox(0,0)[lb]{\smash{{\SetFigFont{6}{7.2}{\rmdefault}{\mddefault}{\updefault}{\color[rgb]{0,0,0}$\lambda_n$}%
}}}}
\put(10801,-3931){\makebox(0,0)[lb]{\smash{{\SetFigFont{6}{7.2}{\rmdefault}{\mddefault}{\updefault}{\color[rgb]{0,0,0}$\hat\lambda_n$}%
}}}}
\end{picture}%

\caption{\textsf{A sample of $n=11$ points.}}
\label{fig:toy_ex}
\end{figure}

We consider a mixture of two uniform univariate densities
$$g_{1,n}(x) = g_1(x) = \1_{[0,1]}(x)\quad\text{and}\quad g_{2,n}(x) =\1_{[1-\lambda_n, 2-\lambda_n]}(x),$$
where we recall that $g_{i,n}$ is the density of the conditional distribution $\mathcal L(X|I=i),i=1,2$. Here $(\lambda_n)_n$ is a non-increasing sequence which tends to $0$ as $n$ goes to infinity. In this parametric situation, a natural way to guess the unobserved label $I_k$ of the observation $X_k$ is to find an estimator $\hat\lambda_n$ of $\lambda_n$ and to predict the labels (see Figure \ref{fig:toy_ex}) according to
\begin{equation}
\label{eq:prev_toyex}
\hat I_k =
  \begin{cases}
    1 &\text{if} \quad X_k \leq 1-\hat\lambda_n\\
    0 &\text{if} \quad 1-\hat\lambda_n < X_k < 1\\
    2 &\text{if} \quad X_k \geq 1.
  \end{cases}
\end{equation}
The accuracy of these predictions depends on the choice of the estimator $\hat\lambda_n$. Here we choose
$\hat\lambda_n=2-X_{(n)}$ where $X_{(n)}=\max_{1\leq k\leq n}X_k$.
Note that in this situation, we have for $i=1,2$
$$\hat I_k=i \Longrightarrow I_k=i,\ \textrm{a.s.}$$
It means that all classified observations (with non-zero estimated label) are well-classified and that misclassified observations are collected in $\cluster{0}$ (see Figure \ref{fig:toy_ex}).

The following proposition establishes a performance bound for the maximal  misclassification error $\varphi_n$ of this clustering procedure.
\begin{pro}
\label{prop:toy_ex}
There exists a  positive constant $A_4$ such that for all $n\geq 1$
  \begin{equation*}
    \varphi_n \leq \lambda_n+A_4\frac{\log n}n.
  \end{equation*}
\end{pro}

Unsurprisingly, $\varphi_n$ decreases as $\lambda_n$ decreases. Moreover, since in most cases of interest, the expected $L_1$-error of $\bfi$ tends to zero much slower than $1/\sqrt{n}$, this property means that, asymptotically, the expected $L_1$-error of $\hfi$ is of the same order as the expected $L_1$-error of $\bfi$ provided $\lambda_n=\mathcal{O}(1/\sqrt{n})$ (see \eqref{eq:borne_er_L1}).

\subsection{A hierarchical clustering algorithm}
\label{sec:supportdisjoint}
Assuming that clusters are defined as connected components of level sets of a density, many authors have studied theoretical properties of various clustering algorithms. For instance, \cite{mahevo09} and \cite{arr11} prove that algorithms based on pairwise distances ($k$-nearest neighbor graph, spectral clustering...) are efficient as soon as these connected components are separated enough. In this section, we extend results of these authors to bound the maximal misclassification error $\varphi_n$ for a hierarchical clustering algorithm.

\subsubsection{The clustering algorithm}
\label{sec:clusmethsuppdisj}
Given $X_1,\hdots,X_n$, we consider a single linkage hierarchical clustering algorithm based on pairwise distances to extract exactly $M$ disjoint clusters $\cluster{1},\ldots,\cluster{M}$ from the observations (see \cite{arr11}). This algorithm consists of finding a data-driven radius $\hat r_n>0$ such that the set
\begin{equation}
\label{eq:unionball}
\bigcup_{k=1}^n B(X_k, \hat r_n)
\end{equation}
has exactly $M$ connected components. Here $B(x,r)$ stands for the closed Euclidean ball with center $x\in\R^d$ and radius $r>0$. Cluster $\cluster{i}$ is then naturally composed by observations $X_k$ which belong to the $i$\textsuperscript{th} connected component of the set~\eqref{eq:unionball}.
\bigskip

The radius $\hat r_n$ can be defined in a formal way to derive statistical properties of the clustering procedure. To this end, we define for each positive real number $r$ the $n\times n$ affinity matrix $A^r=(A^r_{k,\ell})_{1\leq k,\ell\leq n}$ by
\begin{equation}
\label{eq:defmatA}
  A_{k,\ell}^r = \begin{cases}
    1 &\text{if } \norme2{X_k-X_\ell}\leq 2r \iff B(X_k, r)\cap B(X_{\ell},
    r)\neq\emptyset,\\
    0 &\text{otherwise,}
  \end{cases}
\end{equation}
where $\|x\|_2$ stands for the Euclidean norm of $x\in\R^d$. This matrix induces a non-orientated graph on the set $\interval{n}$ and two different observations $X_k$ and $X_{\ell}$ belong to the same cluster if $k$ and $\ell$ belong to the same connected component of the graph. We let $\hat M_r$ be the number of connected components of the graph and we denote by $\cluster{1}(r),\hdots,\cluster{{\hat M_r}}(r)$ the associated
clusters. The radius is selected as follows
$$\hat r_n=\inf\{r>0 : \hat M_r \leq M\}.$$
Note that $\hat r_n$ is well-defined since the random set $\mathcal{R}_M=\{r>0 : \hat M_r \leq M\}$ is lower bounded (by $0$) and non-empty since $r^*=\max_{k,\ell}\norme2{X_k-X_{\ell}}$ always belongs to this set ($\hat M_{r^*}=1$). Moreover, since $r\mapsto\hat M_r$ is non-increasing and right-continuous, one can easily prove that $\hat r_n = \min \mathcal{R}_M$ and $\hat M_{\hat r_n}=M$ almost surely when $n\geq M$. Let $\cluster{1}(\hat r_n),\hdots,\cluster{M}(\hat r_n)$ be the $M$ clusters induced by $A^{\hat r_n}$, the aim is to study the maximal misclassification error \eqref{eq:def_phin} of this clustering algorithm.

\begin{rem}
  The algorithm requires that the connected components of the
  graph induced by the $n\times n$ matrix $A^{r}$ be computed for different values
  of $r$. Some algorithms can be performed to obtain these connected
  components. For instance, we can use the Depth-First search algorithm (see~\cite{coleri90}) which can be performed efficiently in $\gtau (V_n+E_n)$ operations, where $V_n$ and $E_n$ denote
  respectively the number of vertices and edges of the graph.
\end{rem}

\subsubsection{The clustering model}
\label{sec:clust_model}
Recall that the clustering algorithm is performed on the sample $X_1,\hdots,X_n$. To study the maximal misclassification error, some assumptions on the distribution of $X$  are needed.

\paragraph{Assumption 1} Let $g_n$ denotes the probability density of $X$. We assume that there exists a positive sequence $(t_n)_n$ such that the set
\begin{equation}
  \label{eq:ens_niveau}
 \{x\in \R^d:g_{n}(x)\geq t_n\}
\end{equation}
has exactly $M$ disjoint connected compact sets $S_{1,n},\hdots,S_{M,n}$ satisfying, for all $i\in\interval{M}$,
\begin{equation}
\label{eq:cond_numero}
\p(X_1\in\Sin|I_1=i)=\int_{\Sin}g_{i,n}(x)\,\mathrm{d}x >1/2,
\end{equation}
where we recall that $g_{i,n}$ stands for the density of the conditional distribution $\mathcal L(X|I=i),i\in\interval{M}$. We note $S_n=\bigcup_{i=1}^M S_{i,n}$ and
$$\delta_n=\inf_{1\leq i\neq j\leq M}{\textrm{dist}}(S_{i,n},S_{j,n}),$$
where
$$\textrm{dist}(S_{i,n},S_{j,n})=\inf_{x\in S_{i,n}}\inf_{y\in S_{j,n}}\|x-y\|_2.$$

\paragraph{Assumption 2}
\label{hypothèse2}
There exist two positive constants $c_1$ and $c_2$, and a family of $N\in\N^\star$ Euclidean balls $\{B_\ell\}_{\ell=1,\ldots,N}$ with radius $r_n/2$ such that
  \begin{equation*}
    \begin{cases}
      S_n \subset \bigcup_{\ell=1}^N B_\ell\\
      \op{Leb}(S_n) \geq c_1 \sum_{\ell=1}^N \op{Leb}(S_n\cap B_\ell)\\
      \forall\ell=1,\ldots,N,\quad \op{Leb}(S_n\cap B_\ell) \geq c_2 r_n^d,
    \end{cases}
  \end{equation*}
where $\op{Leb}$ denotes the Lebesgue measure on $\R^d$ and $r_n$ is defined by
$$r_n^d = \frac{\tau\log n}{nt_n}\qquad \text{with}\quad\tau>1/c_2.$$
\bigskip

Assumption 1 is classical to study performances of clustering algorithm (see \cite{mahevo09}) or to estimate the number of clusters (see \cite{MR2320821}). It implies that clusters reflect high-density regions separated by low-density regions. Condition \eqref{eq:cond_numero} is required to be sure that the connected components of \eqref{eq:ens_niveau} are correctly indexed. It makes it possible to avoid that most of the observation in $S_{i,n}$ are drawn from $g_{j,n}$ with $j\neq i$.
Assumption 2 is more technical and pertains to the diameter and  regularity of the sets $\Sin$. Our approach consists of identifying sets $\Sin$ with the connected components of $\bigcup_{k=1}^nB(X_k,r)$. Thus, when diameter of $\Sin$ increases, large values of radius $r$ are necessary to connect observations in $\Sin$. However for too large values of $r$, the number of connected components of $\bigcup_{k=1}^nB(X_k,r)$
becomes smaller than $M$ and the method fails. Consequently, we need to constraint the diameter of $\Sin$. This is ensured by assumption 2 since it implies that $S_n$ can be covered by $N$ Euclidean balls such that
\begin{equation}
\label{eq:majorant-N}
N\leq \frac{n}{c_1c_2\tau\log n}.
\end{equation}
Finally, inequality $\op{Leb}(S_n\cap B_\ell) \geq c_2 r_n^d$ in assumption 2 can be seen as a smoothness assumption on the boundaries of $S_n$ (see \cite{MR2463383}).

\begin{rem}
In dimension 1, since each $\Sin$ is connected, it is a segment of the real line. Thus, under assumption 1, its diameter is bounded by $1/t_n$ and assumption 2 is satisfied. For higher dimensions, things turn out to be more complicated. Indeed, even if the measure of the compact set $S_n$ is upper bounded by $1/t_n$, its diameter can be as large as we want.
Consider for example the density $$h_n(x,y)=\1_{[1-1/a_n,a_n]}(x)\1_{[0,1/x^2]}(y),\quad (x,y)\in\R^{+\star}\times\R^+,$$ where $a_n>1$. Since $a_n$ could be chosen to be arbitrarily large, the diameter of $S_n$ could also be arbitrarily large and assumption 2 does not hold. This assumption restricts to some extent the shape of $S_n$. It is satisfied for regular sets such that the diameter does not increase too quickly as $n$ goes to infinity. For example, consider the two dimensional situation where $S_n$ is a rectangle with length $u_n$ and width $v_n$. In such a scenario, one can easily prove that if there exist two positive constants $a_1$ and $a_2$ such that $u_n\geq a_1r_n$ and $v_n\geq a_2 r_n$, then assumption~2 holds. Note also that this assumption is verified for sets $S_n$ that do not depend on the sample size $n$ with smooth boundaries (see~\cite{MR2320821,mahevo09}).
\end{rem}

\begin{figure}[H]
  \centering
\begin{picture}(0,0)%
\includegraphics{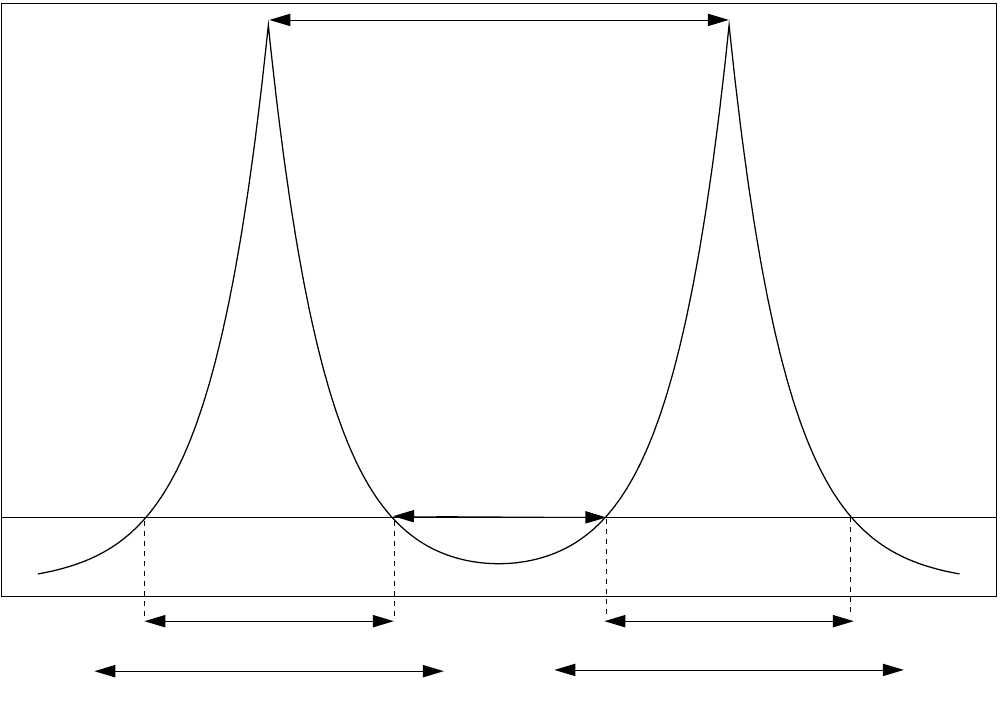}%
\end{picture}%
\setlength{\unitlength}{1579sp}%
\begingroup\makeatletter\ifx\SetFigFont\undefined%
\gdef\SetFigFont#1#2#3#4#5{%
  \reset@font\fontsize{#1}{#2pt}%
  \fontfamily{#3}\fontseries{#4}\fontshape{#5}%
  \selectfont}%
\fi\endgroup%
\begin{picture}(11967,8393)(1273,-8825)
\put(4201,-8161){\makebox(0,0)[lb]{\smash{{\SetFigFont{5}{6.0}{\rmdefault}{\mddefault}{\updefault}{\color[rgb]{0,0,0}$S_{1,n}$}%
}}}}
\put(9901,-8161){\makebox(0,0)[lb]{\smash{{\SetFigFont{5}{6.0}{\rmdefault}{\mddefault}{\updefault}{\color[rgb]{0,0,0}$S_{2,n}$}%
}}}}
\put(3901,-8761){\makebox(0,0)[lb]{\smash{{\SetFigFont{5}{6.0}{\rmdefault}{\mddefault}{\updefault}{\color[rgb]{0,0,0}$S_{1,n}+r_n$}%
}}}}
\put(9601,-8761){\makebox(0,0)[lb]{\smash{{\SetFigFont{5}{6.0}{\rmdefault}{\mddefault}{\updefault}{\color[rgb]{0,0,0}$S_{2,n}+r_n$}%
}}}}
\put(7021,-6346){\makebox(0,0)[lb]{\smash{{\SetFigFont{5}{6.0}{\rmdefault}{\mddefault}{\updefault}{\color[rgb]{0,0,0}$\delta_n$}%
}}}}
\put(6991,-976){\makebox(0,0)[lb]{\smash{{\SetFigFont{5}{6.0}{\rmdefault}{\mddefault}{\updefault}{\color[rgb]{0,0,0}$\ell_n$}%
}}}}
\end{picture}%
  \caption{\textsf{Connected components of level sets for a mixture of Laplace distributions.}}
  \label{fig:dessin_laplace}
\end{figure}

\begin{rem}
\label{exemple:laplace}
Assumption 1 is clearly satisfied when supports of conditional densities $g_{i,n}$ are disjoint. This assumption could also be verified when these supports overlap. As an example, consider the Laplace mixture model:
$$g_{i,n}(x)=\frac1{2\sigma_n}\exp\left(-\frac{|x-\mu_{i,n}|}{\sigma_n}\right),\qquad i=1,2,$$
where $\sigma_n>0$ and $\mu_{i,n}\in\R$ (see Figure \ref{fig:dessin_laplace}). Let $\ell_n=|\mu_{1,n}-\mu_{2,n}|$ be the distance between the two location parameters $\mu_{1,n}$ and $\mu_{2,n}$ and define
$$t_{*,n} = \frac{\sqrt{\alpha_1\alpha_2}}{\sigma_n}\exp\left(-\frac{\ell_n}{2\sigma_n}\right)$$
and
$$t^*_{i,n} = \frac{1}{2\sigma_n}\left(\alpha_i+(1-\alpha_i) \exp\left(-\frac{\ell_n}{\sigma_n}\right)\right),\quad i=1,2.$$
Then direct calculations yield that for any $t_n\in (t_{*,n},t_{1,n}^*\wedge t_{2,n}^*)$, the level set  $\{\alpha_1 g_{1,n}+\alpha_2 g_{2,n}\geq t_n\}$ has exactly $M=2$ connected components provided $\log(\alpha_1/(1-\alpha_1))\in(-\ell_n/\sigma_n, \ell_n/\sigma_n)$.
\end{rem}

\subsubsection{The maximal misclassification error}
The algorithm described in Section~\ref{sec:clusmethsuppdisj} provides a partition of $\{X_1,\hdots,X_n\}$ into $M$ clusters $\cluster{1}(\hat r_n),\hdots,\cluster{M}(\hat r_n)$. To apply Theorem \ref{theo:errL1gen}, we have to find an upper bound of the maximal misclassification error for the predicted rule
$$\hat I_k=i\Longleftrightarrow X_k\in\cluster{i}(\hat r_n).$$
Observe that, for this clustering algorithm, clusters $\cluster{1}(\hat r_n),\hdots,\cluster{M}(\hat r_n)$ defined in Section \ref{sec:clusmethsuppdisj} are arbitrarily indexed. Thus there is no guarantee that the predicted labels are correctly indexed. To circumvent this problem, as suggested in Remark \ref{rem:permutation}, we study the maximal misclassification error up to a permutation of the indexes.
\medskip

The proposed clustering algorithm has been studied by \cite{mahevo09} and \cite{arr11}. They prove that each cluster corresponds to one of the connected components of \eqref{eq:ens_niveau} with high probability in a model similar to ours. In other words, clusters make it possible to identify each connected components of \eqref{eq:ens_niveau}. Even if the identification of these connected components is important in our setting, it is not sufficient since our goal is to find an upper bound of the misclassification error \eqref{eq:phi_perm}. Moreover, since supports of conditional densities $g_{i,n}$ can overlap, observations in the connected components $\Sin$ of \eqref{eq:ens_niveau} are not guaranteed to emerge from the distribution of $\mathcal L(X|I=i)$. This leads us to define
$$\psi_n=\max_{i=1,\hdots,M}\p(X_1\notin (\Sin+r_n)|I_1=i)$$
where for $S\subset\R^d$ and $r>0$
$$S+r=\{x\in\R^d:\exists y\in S\textrm{ such that }\|x-y\|_2\leq r\}.$$
Observe that $\psi_n$ is the maximal probability that an observation from the $i^{th}$ group does not belong to $\Sin+r_n$. This parameter reflects the degree of difficulty for the model to correctly predict  the label of the observations: the larger $\psi_n$, the more difficult it is. We can now set forth the main result of this section.
\begin{theo}
\label{theo:phi_suppdisj}
Suppose that Assumption~1 and Assumption~2 hold. Moreover, if
\begin{equation}
\label{eq:mini_dis_supp}
\delta_n>2r_n =2\left(\frac{\tau\log n}{nt_n}\right)^{1/d},
\end{equation}
then for all $0<a\leq c_2\tau-1$, we have
\begin{equation}
\label{eq:erreur_supp_disj}
\min_{\pi\in\Pi_M}\max_{1,\hdots,n}\p(\pi(\hat I_k)\neq I_k)\leq \frac{A_5}{n^{a}\log n}+(n+2)\psi_n,
\end{equation}
where $A_5$ is positive constant.
\end{theo}
This theorem provides minimal assumptions to make accurate predictions of the labels $I_k$. Inequality~\eqref{eq:mini_dis_supp} gives the minimum distance between the connected components $\Sin$ to make the clustering method efficient. When supports of the conditional densities $g_{i,n}$ are disjoints, it is easily seen that $\psi_n=0$ and $\hat I_k=I_k$ almost surely for $n$ large enough provided inequality \eqref{eq:mini_dis_supp} is satisfied. When the supports overlap, inequality \eqref{eq:erreur_supp_disj} ensures that the algorithm performs well provided the probability $\psi_n$ tends to zero much faster than $1/n$. In the Laplace example presented in Remark~\ref{exemple:laplace}, it can be easily seen that
$$\psi_n=\mathcal O\left(\exp\left(-\frac{\ell_n}{2\sigma_n}\right)\right).$$
It implies that as soon as $\ell_n/\sigma_n\geq3\log(n)/2$, $n\psi(n)\leq n^{-1/2}$ and the kernel density estimates defined in \eqref{eq:def_kern_est} satisfy
$$\min_{\pi\in\Pi_M}\e\|\hat f_{\pi(i)}-f_i\|\leq \e\|\bfi-f_i\|+\frac{A_6}{\sqrt{n}}.$$
\bigskip

Finally, note that when $\psi_n=0$, inequality \eqref{eq:erreur_supp_disj} implies that each cluster $\cluster{i}(\hat r_n)$ belong to one of the connected components of \eqref{eq:ens_niveau} with high probability. This result was obtained by \cite{arr11} in a context similar to ours under assumption \eqref{eq:mini_dis_supp}. Theorem \ref{theo:phi_suppdisj} extends this result for $\psi_n>0$. Note also that proof of this theorem (see Section \ref{sec:proofs}) is different from \cite{arr11} and rely on support density estimation tools proposed by \cite{MR2463383}.

\section{Simulation study}
\label{sec:simu}

In this section, we provide simulation results enlightening the efficiency of the proposed estimator. To this end, $Y$ is simulated from mixtures of univariate Gaussian laws whereas several scenarios on the distribution of $X$ are considered.

\bigskip

To illustrate Theorem \ref{theo:errL1gen} and Theorem  \ref{theo:phi_suppdisj}, we compare the accuracy of our two-step estimate $\hat f_i$ (see \eqref{eq:def_kern_est}) with the accuracy of the oracle estimate $\bar f_i$ (see \eqref{eq:defbarf}). Such comparisons are made in both Sections~\ref{ssec:EM} and \ref{sec:clust_simul}. However, each of these sections focus on special points.
\bigskip

In Section~\ref{ssec:EM}, the two-step estimate is also compared with the classical EM algorithm. Even if this algorithm is known to be efficient under the parametric assumption made on the distribution of $Y$, it does not take advantage of the presence of covariates $X$. It allows our method to outperform the EM algorithm in favorable situations.
\bigskip

In Section~\ref{sec:clust_simul}, different clustering procedures on $X$ are considered on several classical data sets. In particular the behavior of the spectral clustering and the $k$-means algorithm are studied. Both of them are compared with the hierarchical method studied in Section~\ref{sec:supportdisjoint}.

\subsection{Comparison with the EM algorithm}
\label{ssec:EM}
\begin{figure}[H]
    \centering
\begin{picture}(0,0)%
\includegraphics{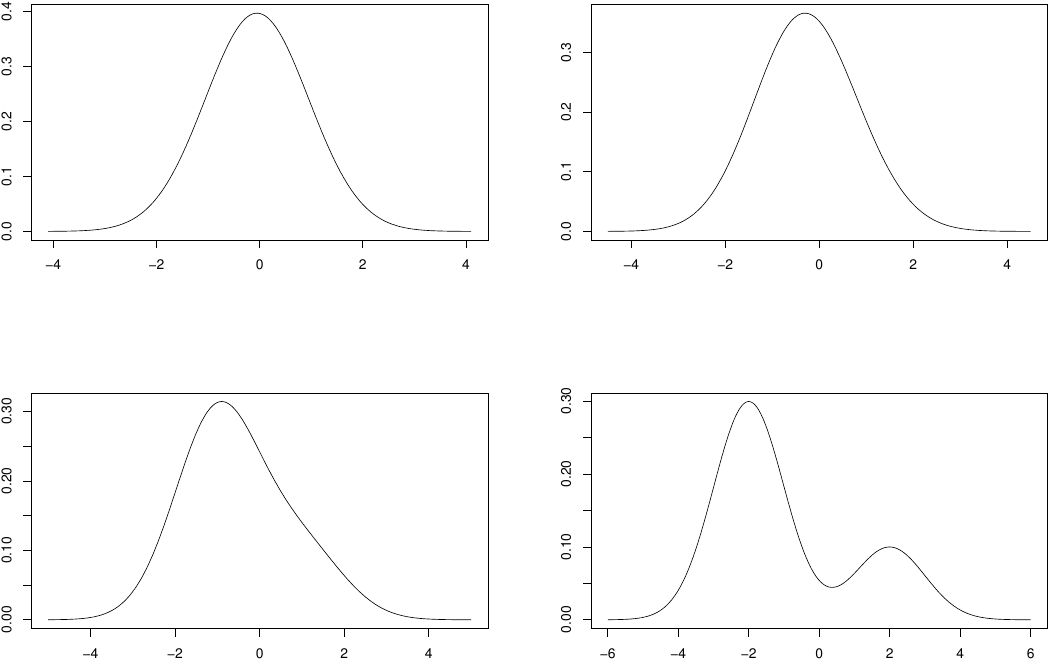}%
\end{picture}%
\setlength{\unitlength}{1579sp}%
\begingroup\makeatletter\ifx\SetFigFont\undefined%
\gdef\SetFigFont#1#2#3#4#5{%
  \reset@font\fontsize{#1}{#2pt}%
  \fontfamily{#3}\fontseries{#4}\fontshape{#5}%
  \selectfont}%
\fi\endgroup%
\begin{picture}(12584,7905)(742,-8139)
\put(5101,-961){\makebox(0,0)[lb]{\smash{{\SetFigFont{6}{7.2}{\rmdefault}{\mddefault}{\updefault}{\color[rgb]{0,0,0}$\Delta=0.1$}%
}}}}
\put(11701,-961){\makebox(0,0)[lb]{\smash{{\SetFigFont{6}{7.2}{\rmdefault}{\mddefault}{\updefault}{\color[rgb]{0,0,0}$\Delta=0.5$}%
}}}}
\put(5101,-5761){\makebox(0,0)[lb]{\smash{{\SetFigFont{6}{7.2}{\rmdefault}{\mddefault}{\updefault}{\color[rgb]{0,0,0}$\Delta=1$}%
}}}}
\put(11701,-5761){\makebox(0,0)[lb]{\smash{{\SetFigFont{6}{7.2}{\rmdefault}{\mddefault}{\updefault}{\color[rgb]{0,0,0}$\Delta=2$}%
}}}}
\end{picture}%

    \caption{{\textsf{Density of $Y$ for various values of $\Delta$.}}}
    \label{fig:densite_Y_simu1}
  \end{figure}

In this simulation section, density of $Y$ is given by
$$f(t)=\frac{3}{4}f_1(t)+\frac{1}{4}f_2(t),\quad t\in\R$$
where $f_1$ and $f_2$ stand for the densities of the normal distribution with mean $-\Delta$ and $\Delta$ and variance $1$. Parameter $\Delta$ measures the separation between the components $f_1$ and $f_2$ (see Figure \ref{fig:densite_Y_simu1}).
Two scenarios are considered for the distribution of $X$. In the first one, conditional densities $g_{i,n},i=1,2$ are uniform univariate densities:
$$g_{1,n}(x)=\mathbb I_{]0,1[}(x)\quad\textrm{and}\quad g_{2,n}(x)=\frac{1}{2}\mathbb I_{]1+\delta_n,3+\delta_n[}(x),\quad x\in\R$$
where $\delta_n>0$ measures the distance between the supports of $g_{1,n}$ and $g_{2,n}$. For the second one, we consider the mixture of Laplace distributions discussed in Section \ref{sec:clust_model}: conditional densities $g_{i,n},i=1,2$ are given by
$$g_{i,n}(x)=\frac{1}{2\sigma_n}\exp\left(-\frac{|x-\mu_{i,n}|}{\sigma_n}\right),\quad i=1,2,$$
where $\sigma_n=1,\mu_{1,n}=1$ and $\mu_{2,n}=\mu_{1,n}+\ell_n$ where $\ell_n>0$. Observe that supports of $g_{i,n}$ are disjoints in the uniform scenario while they overlap in the Laplace example. The separation between these conditional distributions is represented by the location parameters $\delta_n$ and $\ell_n$.
\medskip

For the two proposed scenarios, estimators $\hat f_1$ and $\hat f_2$  defined in \eqref{eq:def_kern_est} are computed using the hierarchical clustering procedure proposed in Section~\ref{sec:supportdisjoint}. These estimates are compared in terms of $L_1$-error with the oracle (but unobservable) estimates $\bar f_1$ and $\bar f_2$ defined in \eqref{eq:defbarf}. Nonparametric kernel estimates $\bar f_i$ and $\hat f_i$ are computed with a Gaussian kernel. Recall that this paper does not put forth any theory for selecting the bandwidth $h$ in an optimal way (see Remark~\ref{rem:bw-choice}). Here we use the default data-driven procedure proposed in the GNU-R library {\bf np} (see \cite{hayrac08}). In addition, these nonparametric density estimates are compared with the EM algorithm (\cite{delaru77}) known to perform well to estimate parameters in a Gaussian mixture model. Formally, we run this algorithm on the sample $Y_1,\hdots,Y_n$ to estimate Gaussian parameters of $f_1$ and $f_2$.  We use the GNU-R library {\bf mclust} and denote by $f_1^{em}$ and $f_2^{em}$ the resulting estimates. They are used as a benchmark. We set $n=300$ and, for the sake of clarity, we present the results regarding $f_1$ only since conclusions are the same for $f_2$. Table \ref{tab:err_L_1_sim1} presents, for different values of $\Delta$, $\delta_n$ and $\ell_n$, the ratio
\begin{equation}
\label{eq:ratio1}
\mathcal R(\tilde f_1)=\frac{\e\|\tilde f_1-f_1\|_1}{\e\|f_1^{em}-f_1\|_1}
\end{equation}
where $\tilde f_1$ is either $\hat f_1$ or $\bar f_1$. Expectations are evaluated over $500$ Monte Carlo replications.
\bigskip

\begin{table}[H]
  \centering
  \begin{tabular}{|c||c|c|c||c|c|c||c|}
\hline
 & \multicolumn{3}{c||}{Uniform:} & \multicolumn{3}{c||}{Laplace:}  &  \multirow{3}{*}{$\mathcal R(\bar f_1)$} \\
 & \multicolumn{3}{c||}{$\mathcal R(\hat f_1)$ for $\delta_n=...$} & \multicolumn{3}{c||}{$\mathcal R(\hat f_1)$ for $\ell_n=...$} &  \\
\cline{2-7}
 & $0.03$ & $0.05$ & $0.1$ & $4.5$ & $5.5$ & $6.5$ &   \\
\hline
\hline
$\Delta=0.1$ & 0.636 & 0.563 & 0.464 & 0.817 & 0.509 & 0.476 & 0.464 \\
$\Delta=0.5$ & 1.156 & 0.923 & 0.679 & 1.261 & 0.749 & 0.692 & 0.679 \\
$\Delta=1$ & 1.772 & 1.288 & 0.844 & 1.769 & 0.954 & 0.869 & 0.843 \\
$\Delta=2$ & 4.243 & 2.876 & 1.702 & 4.298 & 2.093 & 1.830 & 1.701 \\
\hline
  \end{tabular}
  \caption{\textsf{$L_1$-ratio \eqref{eq:ratio1} evaluated over $500$ replications.}}
  \label{tab:err_L_1_sim1}
\end{table}

  \begin{figure}[h]
    \centering
    \includegraphics[width=14cm,height=6cm]{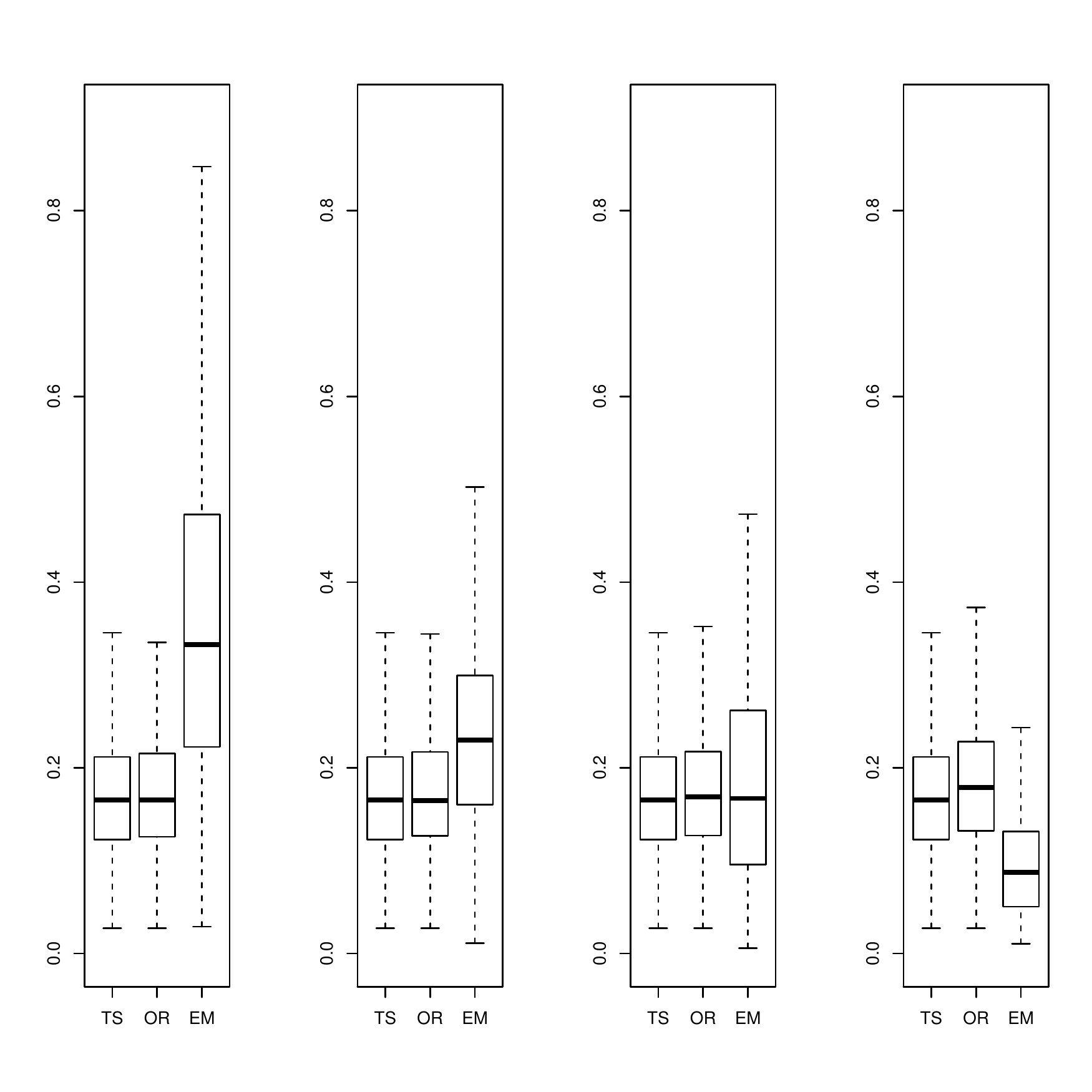}
    \caption{\textsf{Boxplot of the $L_1$-error for the estimate $f_1^{em}$ (EM), the oracle estimate $\bar f_1$ (OR) and the two-step estimate $\hat f_1$ (TS) for the Laplace example. The separation distance $\Delta$ between $f_1$ and $f_2$ vary from 0.1 (left) to 2 (right) and $\ell_n=5.5$.}}
    \label{fig:box_em}
  \end{figure}

\begin{figure}[h]
    \centering
\begin{picture}(0,0)%
\includegraphics{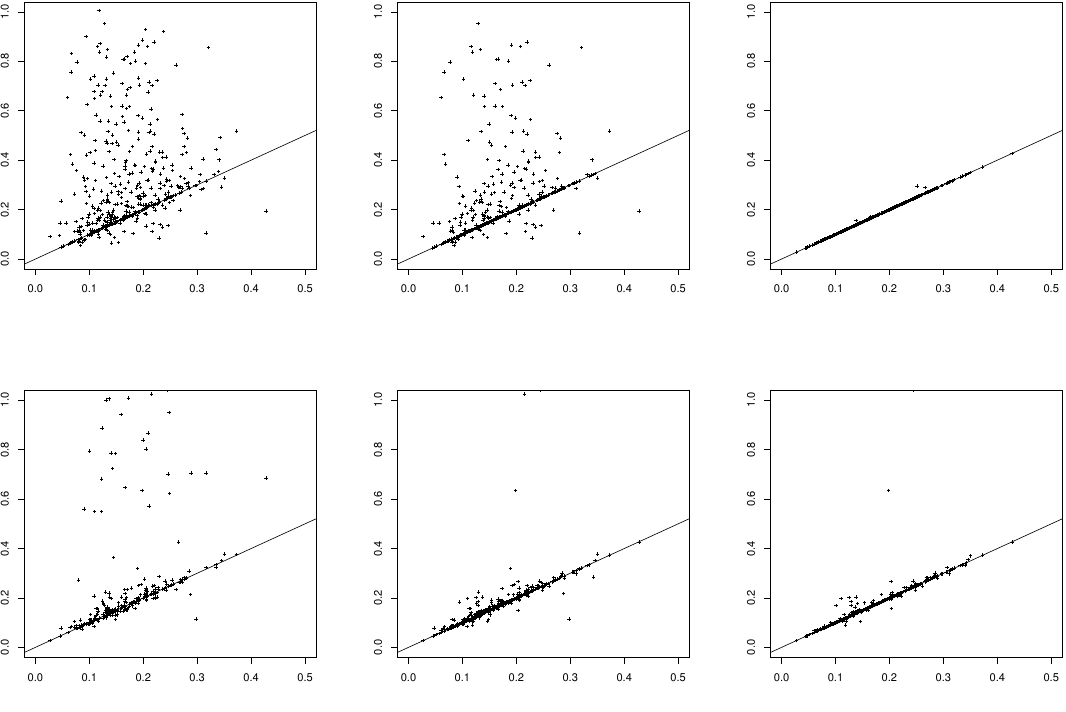}%
\end{picture}%
\setlength{\unitlength}{1579sp}%
\begingroup\makeatletter\ifx\SetFigFont\undefined%
\gdef\SetFigFont#1#2#3#4#5{%
  \reset@font\fontsize{#1}{#2pt}%
  \fontfamily{#3}\fontseries{#4}\fontshape{#5}%
  \selectfont}%
\fi\endgroup%
\begin{picture}(12758,8529)(654,-8627)
\put(11652,-8576){\makebox(0,0)[b]{\smash{{\SetFigFont{5}{6.0}{\sfdefault}{\mddefault}{\updefault}{\color[rgb]{0,0,0}$\ell_n=6.5$}%
}}}}
\put(2698,-3915){\makebox(0,0)[b]{\smash{{\SetFigFont{5}{6.0}{\sfdefault}{\mddefault}{\updefault}{\color[rgb]{0,0,0}$\delta_n=0.03$}%
}}}}
\put(7175,-3915){\makebox(0,0)[b]{\smash{{\SetFigFont{5}{6.0}{\sfdefault}{\mddefault}{\updefault}{\color[rgb]{0,0,0}$\delta_n=0.05$}%
}}}}
\put(11652,-3915){\makebox(0,0)[b]{\smash{{\SetFigFont{5}{6.0}{\sfdefault}{\mddefault}{\updefault}{\color[rgb]{0,0,0}$\delta_n=0.1$}%
}}}}
\put(2698,-8576){\makebox(0,0)[b]{\smash{{\SetFigFont{5}{6.0}{\sfdefault}{\mddefault}{\updefault}{\color[rgb]{0,0,0}$\ell_n=4.5$}%
}}}}
\put(7175,-8576){\makebox(0,0)[b]{\smash{{\SetFigFont{5}{6.0}{\sfdefault}{\mddefault}{\updefault}{\color[rgb]{0,0,0}$\ell_n=5.5$}%
}}}}
\end{picture}%
    \caption{{\textsf{$L_1$-error of $\bar f_1$ (x-axis) and $\hat f_1$ (y-axis) for the uniform (up) and Laplace (down) example.}}}
    \label{fig:scatter}
  \end{figure}

As expected, the performances of the EM algorithm clearly depend on the separation distance between the target densities $f_1$ and $f_2$. For large  $\Delta$ values, parametric estimates resulting from the EM algorithm outperform the nonparametric estimates proposed in this paper (e.g. $\Delta=2$ in Figure~\ref{fig:box_em}). This is not the case when $f_1$ is closed to $f_2$: $L_1$-performance of $\hat f_1$ over $f_1^{em}$ is significantly better for $\Delta=0.1$ and $\Delta=0.5$ and roughly similar for $\Delta=1$. Note also that the $L_1$-error of $\hat f_1$ does not depend on $\Delta$ (see Figure \ref{fig:box_em}). Figure \ref{fig:scatter} displays scatterplots of the $L_1$-error of $\hat f_1$ versus those of the oracle $\bar f_1$ for $\Delta=1$. As proved in Theorem \ref{theo:errL1gen}, most points are above the diagonal.  The distance from a point to the first bisector measures to some extent the distance between $\hat f_1$ and $\bar f_1$ in terms of $L_1$-error. The closer to the bisector, the better $\hat f_1$.  In other words, this distance represents the performance of the clustering algorithm. We observe that points move closer to the first bisector as separation parameters $\delta_n$ and $\ell_n$ increase.
As explained in Theorem \ref{theo:phi_suppdisj}, performances of the hierarchical clustering algorithm depend on the separation parameters $\delta_n$ and $\ell_n$: when these parameters increase, performances of $\hat f_1$ become similar to those of the oracle $\bar f_1$. Indeed, in our simulations, we observe that $L_1$-error of $\hat f_1$ and $\bar f_1$ are quite the same for $\delta_n=0.1$ (resp. $\ell_n=6.5$) in the uniform case (resp. Laplace case).

\subsection{A comparison of clustering algorithms}
\label{sec:clust_simul}
As discussed in Section \ref{sec:cluster}, any clustering algorithm could be applied in practice. However, it is clear that $L_1$-performances of the proposed estimate depend largely on the performances of the clustering method.
The problem is to find the appropriate clustering algorithm according to the covariates $X$. In this section, we propose to compare three standard clustering procedures: the hierarchical clustering algorithm presented in Section \ref{sec:supportdisjoint}, the spectral clustering algorithm performed with a Gaussian kernel (see \cite{arr11}) and the $k$-means algorithm.
\medskip

\begin{figure}[H]
  \centering
  \includegraphics[width=14cm,height=6cm]{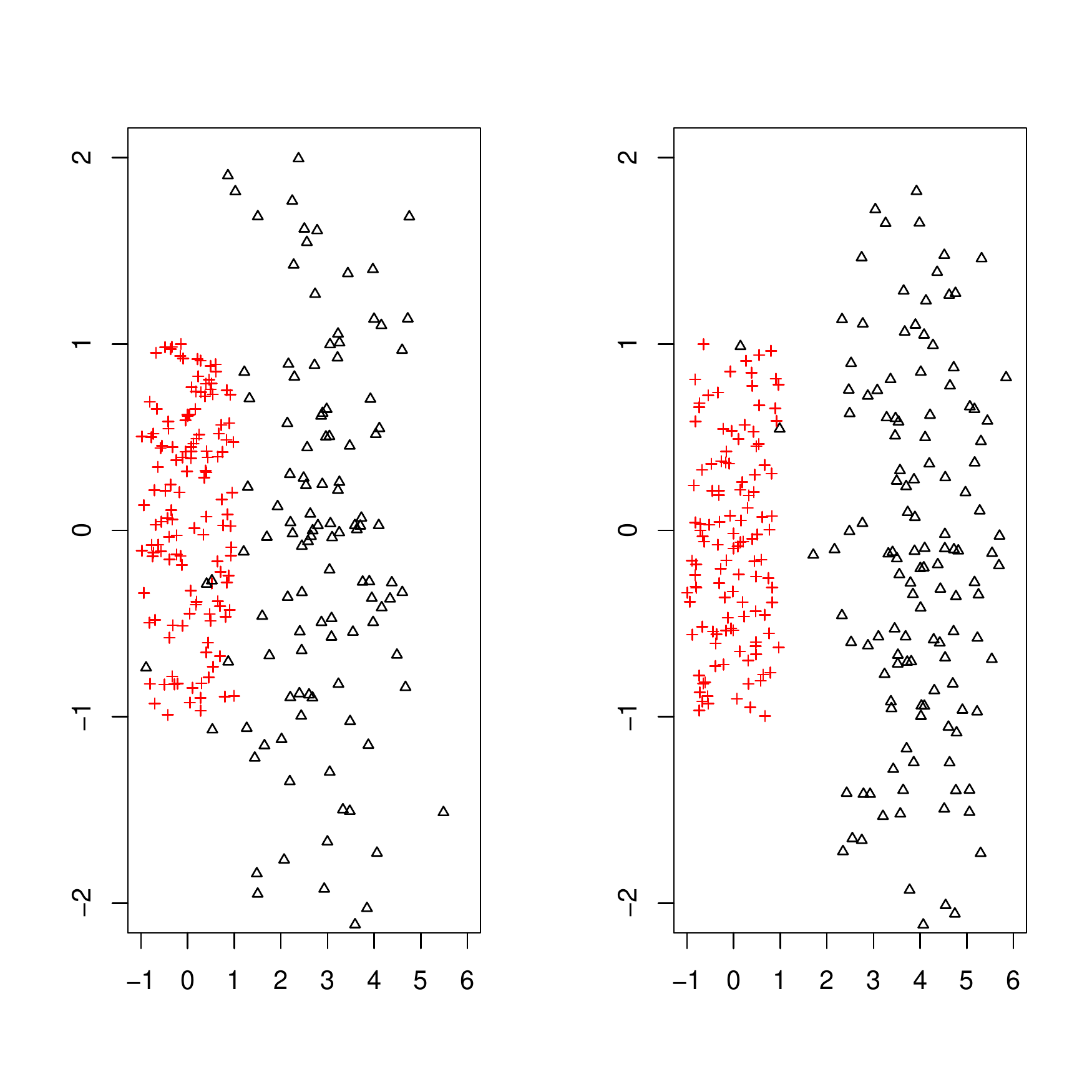}
  \caption{\textsf{A sample of $n=250$ observations for the ``Circle-Square'' model with $a=3$ (left) and $a=4$ (right).}}
  \label{fig:sample_circle_square}
\end{figure}

\begin{figure}[H]
  \centering
  \includegraphics[width=14cm,height=6cm]{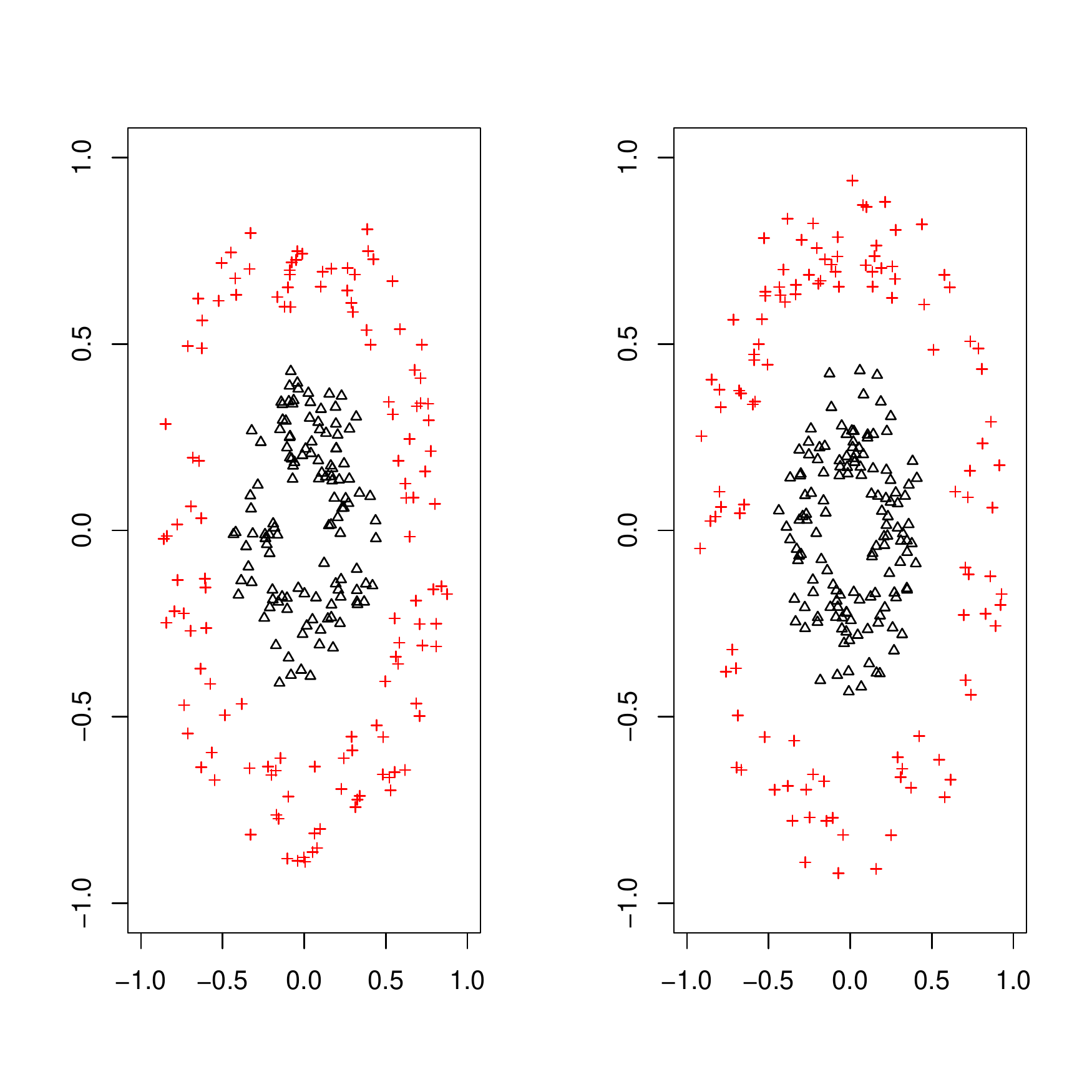}
  \caption{\textsf{A sample of $n=250$ observations for the ``Concentric circles'' model with $r_2=0.75$ (left) and $r_2=0.80$ (right).}}
  \label{fig:sample_concen_circles}
\end{figure}

The model is as follows. The density of $Y$ is now given by
$$f(t)=\frac{1}{2}f_1(t)+\frac{1}{2}f_2(t),\quad t\in\R$$
where $f_1$ and $f_2$ stand for the densities of the normal distribution with mean $-1$ and $1$ and variance $1$.  Here, random variable $X$ takes values in $\R^2$ and we  again consider two scenarios for its distribution:
\begin{itemize}
\item ``Circle-Square'' model (see \cite{bau09}): $g_{1,n}$ is the density of the Gaussian distribution with mean $(a,0)$ and identity variance covariance matrix; $g_{2,n}$ is the density of the uniform distribution over the square $[-1,1]^2$ (see Figure \ref{fig:sample_circle_square}).
\item ``Concentric circles'' model (see \cite{ngjowe02}): $g_{1,n}$ is the density of the uniform distribution over $\mathcal C(0,r_1+\varepsilon,r_1-\varepsilon)$ and $g_{2,n}$ represents the uniform distribution over $\mathcal C(0,r_2+\varepsilon,r_2-\varepsilon)$, where for $r>0$ and $\varepsilon>0$ $\mathcal C(0,r+\varepsilon,r-\varepsilon)$ represents the set between circles with center 0 and radius $r+\varepsilon$ and $r-\varepsilon$ (see Figure \ref{fig:sample_concen_circles}). We fix $r_1=0.3$, $\varepsilon=0.15$ and consider many values for $r_2$ such that $r_2>r_1+2\varepsilon$.
\end{itemize}
The difficulty encountered in identifying each group depends on parameters $a$ and $r_2$. The smaller $a$ and $r_2$, the harder to identify the clusters.

\begin{figure}[H]
    \centering
    \includegraphics[width=15cm,height=7cm]{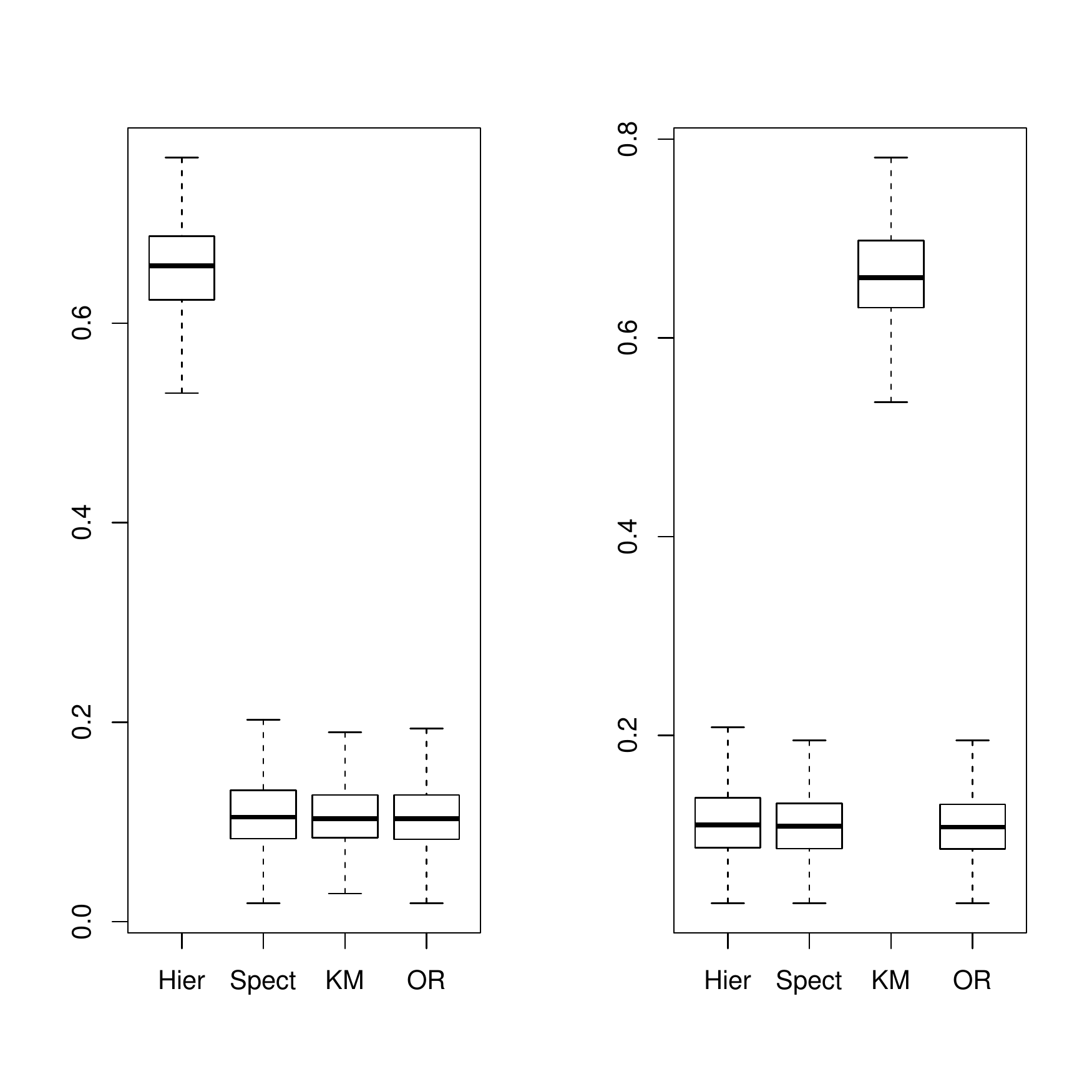}
    \caption{\textsf{Boxplot of the $L_1$-error for the oracle estimate $\bar f_1$ (OR) and two-step estimator $\hat f_1$ using the hierarchical algorithm (Hier), spectral clustering algorithm (Spect) and $k$-means algorithm (KM). Results are for ``Circle-Square'' dataset with $a=4$ and $n=500$ (left) and ``Concentric circles'' dataset with $r_2=0.75$ and $n=500$ (right).}}
    \label{fig:box_clust}
  \end{figure}

For the two described examples, we use the two-step kernel density estimator for three clustering algorithms: hierarchical, spectral and $k$-means. The resulting estimates are compared with the oracle estimates $\bar f_1$ and $\bar f_2$. We keep the same setting as above to compute estimates $\hat f_1$ and $\hat f_2$: Gaussian kernel and bandwidth selected with the library {\bf np}. For the sake of clarity, we again present only results on $\hat f_1$ since we observe the same conclusions for $\hat f_2$. Table \ref{tab:ex:circlequare} and Table \ref{tab:ex:concentcircle} present the ratio
\begin{equation}
  \label{eq:ratio2}
\mathcal R(\hat f_1)=\frac{\e\|\hat f_1-f_1\|_1}{\e\|\bar f_1-f_1\|_1},
\end{equation}
for many values of $a$, $r_2$ and $n$. Expectations are evaluated over $500$ Monte-Carlo replications and Figure \ref{fig:box_clust} presents boxplots of the $L_1$-error of the different estimates. For each replications, we also compute the error of the clustering procedure
$$\frac{1}{n}\sum_{k=1}^n \1_{\hat I_k\neq I_k}$$
and we display in Table \ref{tab:ex:circlequare} and Table \ref{tab:ex:concentcircle} this error term averaged over the 500 replications (it is denoted $err_n$). Observe that this term is closely related to the maximal misclassification error $\varphi_n$.

\begin{table}[H]
  \centering
  \begin{tabular}{|c|c||c|c|c|c|c|c|}
\hline
\multicolumn{2}{|c||}{}  & \multicolumn{2}{c|}{Hier.} & \multicolumn{2}{c|}{Spect.} & \multicolumn{2}{c|}{k-means} \\
\multicolumn{2}{|c||}{} & \multicolumn{1}{c}{$\mathcal R(\hat f_1)$} & $err_n$ & \multicolumn{1}{c}{$\mathcal R(\hat f_1)$} & $err_n$ & \multicolumn{1}{c}{$\mathcal R(\hat f_1)$} & $err_n$ \\
\hline\hline
\multirow{2}{*}{$a=3$} & $n=250$ & 4.680 & 0.475 & 1.748 & 0.121 & 1.047 & 0.043  \\
& $n=500$ & 6.370 & 0.483 & 2.265 & 0.126 & 1.034 & 0.043 \\
\hline
\multirow{2}{*}{$a=4$} & $n=250$ & 3.565 & 0.382 & 1.107 & 0.018 & 1.005 & 0.013 \\
& $n=500$ & 5.688 & 0.449 & 1.190 & 0.023 & 1.000 & 0.013 \\
\hline
\multirow{2}{*}{$a=5$} & $n=250$ & 1.285 & 0.067 & 0.999 & 0.001 & 0.997 & 0.003 \\
& $n=500$ & 1.897 & 0.130 & 0.999 & 0.001 & 1.000 & 0.003 \\
\hline
  \end{tabular}
  \caption{\textsf{Error ratio \eqref{eq:ratio2} evaluated over 500 Monte Carlo replications for the ``Circle-Square'' example.}}
  \label{tab:ex:circlequare}
\end{table}

\begin{table}[H]
  \centering
  \begin{tabular}{|c|c||c|c|c|c|c|c|}
\hline
\multicolumn{2}{|c||}{}  & \multicolumn{2}{c|}{Hier.} & \multicolumn{2}{c|}{Spect.} & \multicolumn{2}{c|}{k-means} \\
\multicolumn{2}{|c||}{} & \multicolumn{1}{c}{$\mathcal R(\hat f_1)$} & $err_n$ & \multicolumn{1}{c}{$\mathcal R(\hat f_1)$} & $err_n$ & \multicolumn{1}{c}{$\mathcal R(\hat f_1)$} & $err_n$ \\
\hline\hline
\multirow{2}{*}{$r_2=0.75$} & $n=250$ & 4.040 & 0.349 & 2.776 & 0.195 & 4.568 & 0.468  \\
& $n=500$ & 1.197 & 0.021 & 1.013 & 0.001 & 5.993 & 0.478 \\
\hline
\multirow{2}{*}{$r_2=0.80$} & $n=250$ & 1.852 & 0.105 & 1.433 & 0.049 & 4.556 & 0.467 \\
& $n=500$ & 1.010 & 0.001 & 1.000 & 0.000 & 5.986 & 0.477 \\
\hline
  \end{tabular}
  \caption{\textsf{Error ratio \eqref{eq:ratio2} evaluated over 500 Monte Carlo replications for the ``Concentric circles'' example.}}
  \label{tab:ex:concentcircle}
\end{table}
As proved in Theorem \ref{theo:errL1gen}, performances of $\hat f_1$ depend on the accuracy of the clustering approach: the lower $err_n$, the better $\hat f_1$. For the ``Circle Square'' dataset, unsurprisingly $k-$means algorithm overperforms the two other clustering methods. Indeed, $k$-means is well appropriate to this dataset since clusters can be identified by their distances to two particular points (the centers of the uniform and Gaussian distributions). It is not the case for the ``Concentric circle'' dataset where estimates defined from hierarchical and spectral clustering algorithms achieve the best estimated $L_1$-error.

\section{Application to electricity distribution}
\label{sec:erdf}

\subsection{Context of the study}

ERDF is the contract-holder of the public electricity distribution network in France. ERDF is in charge of operating, maintaining and developing the network. With 36,000 employees and 35 million customers served over 34,220 communes, ERDF is the largest electricity distributor in Europe. It operates more than 1.3 million km in power lines and runs more than 11 million operations per year. ERDF also plays an essential role in ensuring the proper functioning of the competitive electricity market by providing quality electricity supply among the best in Europe, and serving all network users without resorting to guaranteeing discriminatory practices.

\bigskip

In recent years, the electricity sector  has entered a period of profound changes resulting from the emergence of decentralized and intermittent (wind, solar) means of production and  new electricity uses (e.g.\@ electric vehicle). The increasing integration of these new means of production and new uses has a major impact on ERDF's core business: connecting new users (producers, terminals electric vehicle), and adaptating rules of conduct and network planning/investment to meet the new specifications. ERDF has initiated its digital transformation plan so as to take advantage of new information technologies, and by meeting its new challenges, offer better public service.
\bigskip

ERDF launched the ``smart grid'' experimental programs in order to run the network with more flexibility and efficiency. To do so, these programs use detailed network status and mine/produce information from different users. These more detailed data (including from a new generation of electricity meters, called smart meters) will accordingly be used to improve network monitoring (predictive maintenance).
\bigskip

In this section we focus on the detection of customers who experience a significant decrease in consumption, for a given period of time, {\it i.e.\@}, a period when overall malfunction of the network could be observed. This will make it possible to better understand the origin of dysfunctions and thus better forecast network operation. For this study, we have the benefit of a set of consumption curves for 226 customers with observations taken at regularly spaced instants. Based on the observation of the individual consumption curves, we can cluster individuals into two groups (those who have suffered an abnormal decline and the others) and estimate, in each group, distributions of many variables using the approach proposed in this paper.

\subsection{Application of the two-step estimator}

The consumption curves of $n=226$ ERDF's customers are observed at 9 regularly spaced instants $t_1,\ldots,t_9$. The time interval $[t_1, t_9]$ covers a known period of disruption between times $t_4$ and $t_6$. The observations consist of $n$ vectors $\bZ_k=(Z_{k1},\hdots,Z_{k9})\in\R^9$ where $Z_{kj}$ stands for the consumption of user $k$ at time $t_j$.
\bigskip

Since ERDF is interested in comparing the behavior of customers of both sub-populations (those who have suffered from the disruption and others) before and after the disruption period, we consider $6$ different  variables $Y^{(j)}$ in relation with the consumption around the disruption period. These variables, presented below, are observed for each customer and thus are defined for any $k\in\interval{n}$.

\begin{figure}[h]
  \centering
  \includegraphics[width=10cm,height=6cm]{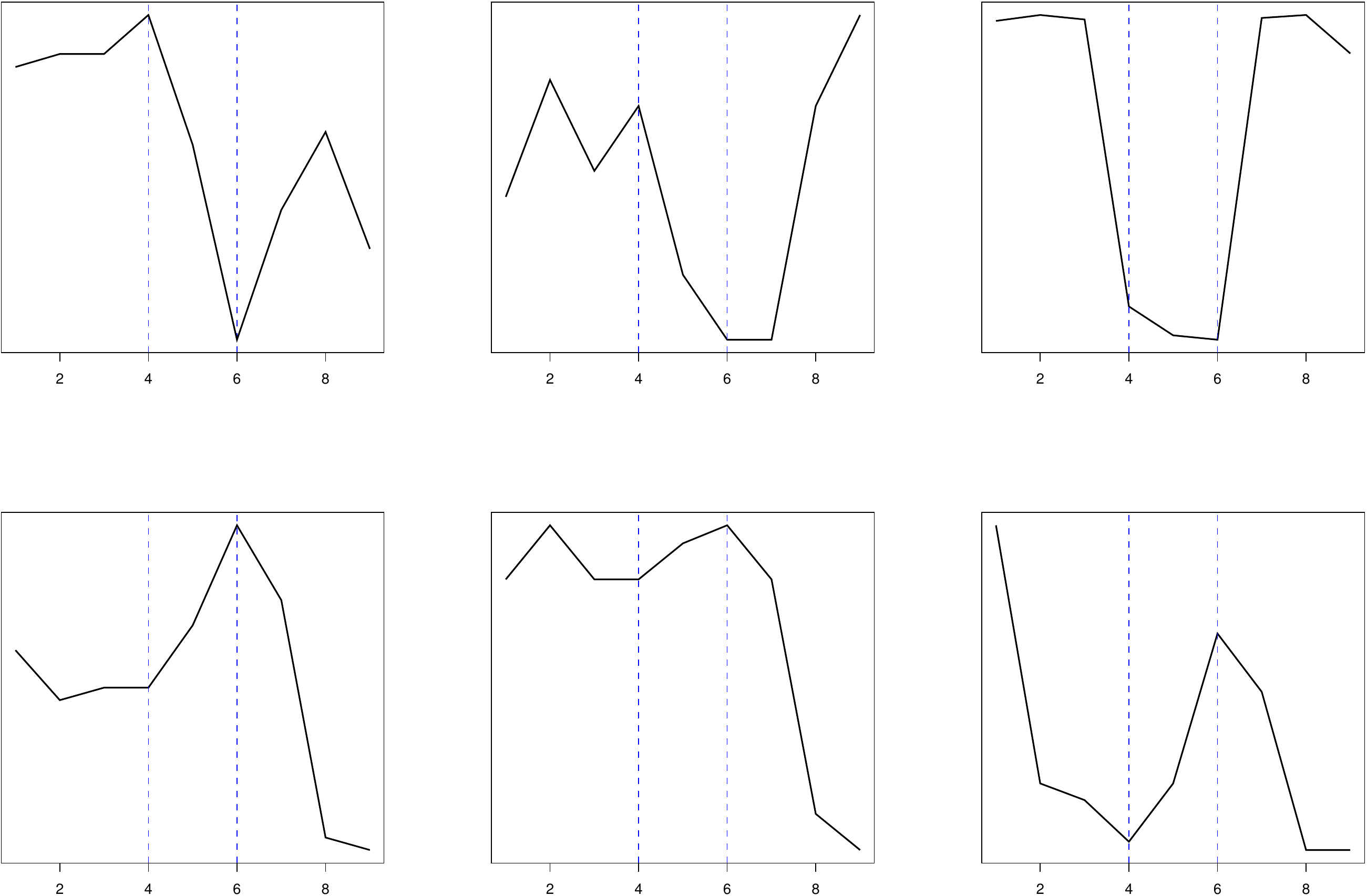}
  \caption{\textsf{Consumptions of users suspected to be affected (up) or not (down) by the perturbation.}}
  \label{fig:courbes_ERDF}
\end{figure}

\begin{enumerate}
\item Average consumptions before, during  and after the disruption period defined by:
$$Y_k^{(1)}=\frac{Z_{k1}+Z_{k2}+Z_{k3}}{3},\quad Y_k^{(2)}=\frac{Z_{k4}+Z_{k5}+Z_{k6}}{3}$$
$$\textrm{and}\quad Y_k^{(3)}=\frac{Z_{k7}+Z_{k8}+Z_{k9}}{3};$$
\item Evolutions of consumption around the disruption period defined by:
$$Y_k^{(4)}=\frac{Y_k^{(2)}-Y_k^{(1)}}{Y_k^{(1)}},\quad Y_k^{(5)}=\frac{Y_k^{(3)}-Y_k^{(1)}}{Y_k^{(1)}}\quad\textrm{and}\quad Y_k^{(6)}=\frac{Y_k^{(3)}-Y_k^{(2)}}{Y_k^{(2)}}.$$
\end{enumerate}
Let $I$ be the random variable taking value 1 if a customer has been affected by the disruption, 2 otherwise. If we denote by $f_1^{(j)}$ and $f_2^{(j)}$ the conditional densities of $\mathcal L(Y^{(j)}|I=1)$ and $\mathcal L(Y^{(j)}|I=2)$, the problem is to compare $f_1^{(j)}$ with $f_2^{(j)}$ for each $j\in\interval{6}$. Even if ERDF can measure consumptions during the disruption period (between $t_4$ and $t_6$), it does not have the capacity to identify consumers affected by the perturbation. It means that random variables $I_k,k=1,\hdots,n$ are not observed. However, we know that users impacted by the disruption posted a decline in consumption during $t_4$ and $t_6$. Figure \ref{fig:courbes_ERDF} provides examples of customers potentially affected by the disruption (for confidentiality reasons, representations are anonymous and scales of power are not specified).

Using the approach developed in this paper, we first have to identify users impacted by the disruption with a clustering algorithm.
As the disruption influences the consumptions of user $k$ between $t_4$ and $t_6$ we define $X_k=(X_{k1},X_{k2}),k=1,\hdots,n$ with
$$X_{k1}=\min\left(v_{k,54},v_{k,65}\right),\quad X_{k2}=v_{k,54}+v_{k,65}$$
where
$$v_{k,ij}=\frac{Z_{kj}-Z_{ki}}{Z_{ki}}, \quad 1\leq i,j\leq 9.$$
Observe that $v_{k,ij}$ measures the relative variation of consumption for user $k$ between $t_i$ and $t_j$. It follows that $X_k=(X_{k1},X_{k2})$ captures the development of consumption  of user $k$ during the disruption period. We use these covariates to cluster users into two groups: the first contains consumers assumed to be affected by the disruption, the second contains the others.

\begin{figure}[h]
  \centering
\begin{picture}(0,0)%
\includegraphics{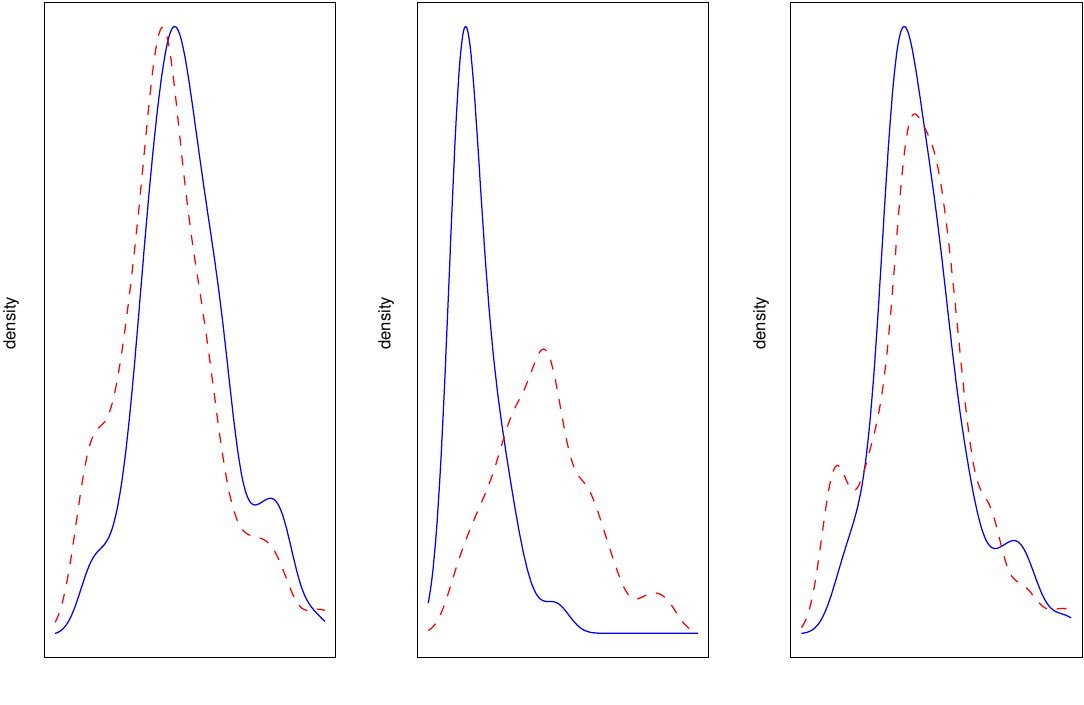}%
\end{picture}%
\setlength{\unitlength}{1579sp}%
\begingroup\makeatletter\ifx\SetFigFont\undefined%
\gdef\SetFigFont#1#2#3#4#5{%
  \reset@font\fontsize{#1}{#2pt}%
  \fontfamily{#3}\fontseries{#4}\fontshape{#5}%
  \selectfont}%
\fi\endgroup%
\begin{picture}(12995,8441)(417,-8539)
\put(6901,-8461){\makebox(0,0)[lb]{\smash{{\SetFigFont{5}{6.0}{\rmdefault}{\mddefault}{\updefault}{\color[rgb]{0,0,0}$Y_1^{(2)}$}%
}}}}
\put(11401,-8461){\makebox(0,0)[lb]{\smash{{\SetFigFont{5}{6.0}{\rmdefault}{\mddefault}{\updefault}{\color[rgb]{0,0,0}$Y_1^{(3)}$}%
}}}}
\put(2401,-8461){\makebox(0,0)[lb]{\smash{{\SetFigFont{5}{6.0}{\rmdefault}{\mddefault}{\updefault}{\color[rgb]{0,0,0}$Y_1^{(1)}$}%
}}}}
\end{picture}%

  \caption{\textsf{Kernel estimates $\hat f_1^{(j)}$ (solid lines) and $\hat f_2^{(j)}$ (dashed lines) for $j=1$ (left), 2 (center) and 3 (right).}}
  \label{fig:res_erdf1}
\end{figure}

\begin{figure}[H]
  \centering
\begin{picture}(0,0)%
\includegraphics{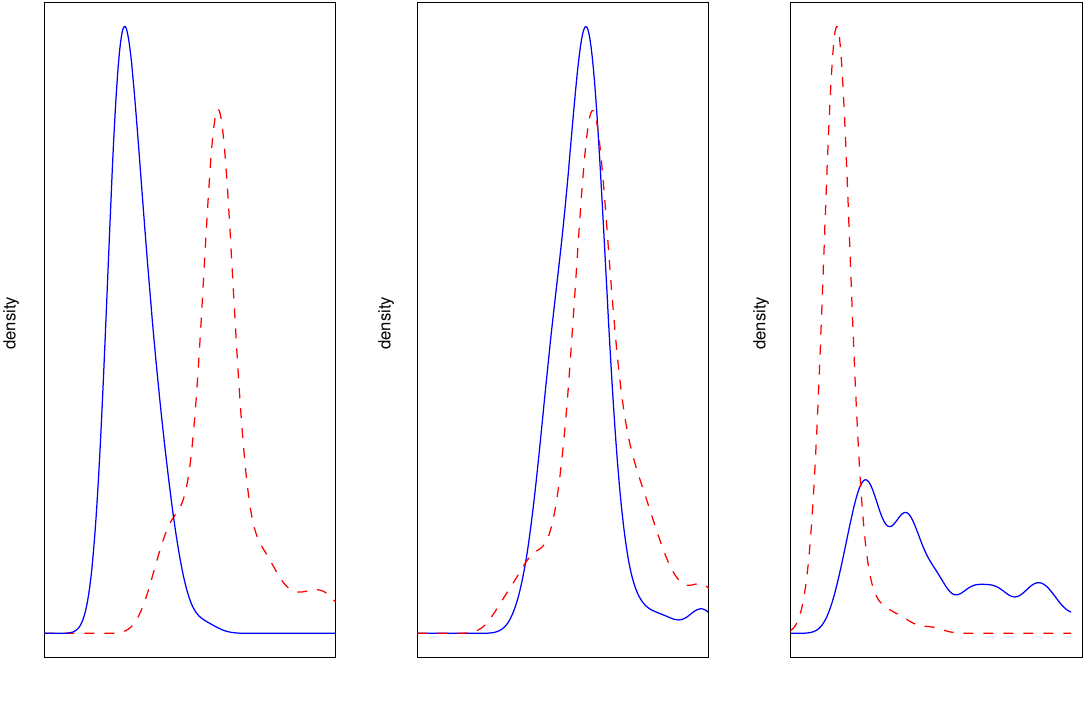}%
\end{picture}%
\setlength{\unitlength}{1579sp}%
\begingroup\makeatletter\ifx\SetFigFont\undefined%
\gdef\SetFigFont#1#2#3#4#5{%
  \reset@font\fontsize{#1}{#2pt}%
  \fontfamily{#3}\fontseries{#4}\fontshape{#5}%
  \selectfont}%
\fi\endgroup%
\begin{picture}(12995,8441)(417,-8539)
\put(2401,-8461){\makebox(0,0)[lb]{\smash{{\SetFigFont{5}{6.0}{\rmdefault}{\mddefault}{\updefault}{\color[rgb]{0,0,0}$Y_1^{(4)}$}%
}}}}
\put(6901,-8461){\makebox(0,0)[lb]{\smash{{\SetFigFont{5}{6.0}{\rmdefault}{\mddefault}{\updefault}{\color[rgb]{0,0,0}$Y_1^{(5)}$}%
}}}}
\put(11401,-8461){\makebox(0,0)[lb]{\smash{{\SetFigFont{5}{6.0}{\rmdefault}{\mddefault}{\updefault}{\color[rgb]{0,0,0}$Y_1^{(6)}$}%
}}}}
\end{picture}%
  \caption{\textsf{Kernel estimates $\hat f_1^{(j)}$ (solid lines) and $\hat f_2^{(j)}$ (dashed lines) for $j=4$ (left), 5 (center) and 6 (right).}}
  \label{fig:res_erdf2}
\end{figure}

\bigskip
Two clustering algorithms have been tested: the hierarchical method studied in section \ref{sec:supportdisjoint} and the $k$-means algorithm. Since these methods lead to approximately the same clusters, we only present results for the hierarchical method. Figures \ref{fig:res_erdf1} and \ref{fig:res_erdf2} present kernel density estimates \eqref{eq:def_kern_est} of conditional densities $f_1^{(j)}$ and $f_2^{(j)}$ for $j\in\interval{6}$. Parameters (bandwidth and kernel) of the kernel estimates are chosen as in the simulation part. For confidentiality reasons, scales of power are again not specified.
\bigskip

Figure \ref{fig:res_erdf1} strongly supports the idea that the clustering procedure allows to correctly identify users impacted by the disruption. Indeed we observe that the average consumption during the disruption period is lower for consumers in the first group (second graph in Figure \ref{fig:res_erdf1}). We can also observe that average consumptions are quite the same for the two groups before and after the disruption period. It means that users impacted by the perturbation do not  over-consume after the disruption period. This conclusion is also supported by the second graph in Figure~\ref{fig:res_erdf2}: distributions representing the evolution of consumptions are similar for the two clusters.

\section{Conclusion}
\label{sec:conclusion}
This paper provides a new framework to estimate conditional densities in mixture models in the presence of covariates. To our knowledge, no clear probabilistic model has been proposed to take into account of the presence of covariates. The model we consider includes such covariates and Theorem \ref{theo:errL1gen} precisely describes the interest of a preliminary clustering step on these covariates to estimate components of the mixture model. It is shown that the performances of these estimates depend on the maximal misclassification error  \eqref{eq:def_phin} of the clustering algorithm. This criterion is natural to measure performances of clustering algorithms but, as far as we know, it has not been addressed before. We obtain non-asymptotic upper bounds of this error term is section \ref{sec:supportdisjoint} for a particular hierarchical algorithm. This algorithm is not new but it has not been studied in this context. Results obtained for this algorithm could be extended to other clustering algorithms based on pairwise distances such as spectral clustering (\cite{arr11}) or on clustering methods based on neighborhoods graphs (\cite{mahevo09}). Even if main contributions of this work are theoretical, both the simulation study and the application on real data enlighten the efficiency of the proposed estimator in the presence of covariates.

\paragraph{Acknowledgement}
We would like to thank the editor, an associate editor as well as two anonymous referees for very thoughtful and detailed comments. We are also grateful to Datastorm and ERDF for providing us with the data-set used in the case study.

\section{Proofs}
\label{sec:proofs}
\subsection{Proof of Theorem~\ref{theo:errL1gen}}
We first prove inequality~\eqref{eq:borne_er_L1}. Since
\begin{align*}
\e\bnorme1{\hfi-\ffi} &\leq
 \e\bnorme1{\bfi-\ffi}+\e\bnorme1{\hfi-\bfi},
\end{align*}
we need only find an upper bound of the second term in the right-hand side  of the previous inequality.
Since $\bfi=0$ when $N_i=0$ and $\pnorme1{\hfi}=\pnorme1K$, we have
\begin{align*}
\e\bnorme1{\hfi-\bfi}
&\leq \e\left(\bnorme1{\hfi}\1_{N_i=0}\right)+\e\bnorme1{(\hfi-\bfi)\1_{N_i>0}}\\
& \leq \pnorme1K(1-\alpha_i)^n+\e\bnorme1{(\hfi-\bfi)\1_{N_i>0}}.
\end{align*}
For the sake of readability, let $\te$ denote the conditional expectation with respect to $(I_1,\hdots,I_n)$ and $\tte$ the conditional expectation with respect to $(I_1,\hdots,I_n,X_1,\hdots,X_n)$. Moreover, let
\begin{align*}
  A_i(t) &= \big(\hfi(t)-\bfi(t)\big)\1_{Ni>0}\\
  &=\sum_{k=1}^nK_{h}(t,Y_k)
  \bigg(\frac{\1_{\{i\}}(\hat I_k)}{\hat
      N_i}-\frac{\1_{\{i\}}(I_k)}{N_i}\bigg)\1_{N_i>0}.
\end{align*}
Using these notations it is easily seen that
\begin{equation}\label{eq:proof-decomp}
  \e\norme1{(\hfi-\bfi)\1_{N_i>0}} = \e\te \int_\R \tte|A_i(t)| dt.
\end{equation}

Since, for all $y\in\R$ we have $\int_\R |K_h(t,y)|dt=\pnorme1K$, we deduce that
\begin{align*}
\int_\R\tte|A_i(t)|\, dt & \leq \sum_{k=1}^n\tte\left(\int_\R|K_{h}(t,Y_k)|\,dt\right)\left|\frac{\1_{\{i\}}(\hat I_k)}{\hNi}-\frac{\1_{\{i\}}(I_k)}{N_i}\right|\\
& \leq \norme1{K}\sum_{k=1}^n\left|\frac{\1_{\{i\}}(\hat I_k)}{\hNi}-\frac{\1_{\{i\}}(I_k)}{N_i}\right|.
 \end{align*}
Thus
\begin{equation}\label{eq:proof-1}
\te\int_\R\tte|A_i(t)|\,dt\leq
\frac{\norme1{K}}{N_i}\te\left[\frac{1}{\hNi}\sum_{k=1}^n|N_i\1_{\{i\}}(\hat
  I_k)-\hNi\1_{\{i\}}(I_k)|\right].
\end{equation}
Moreover, inserting $\hNi\1_{\{i\}}(\hat I_k)$ in the previous expectation, we obtain
\begin{align}
\te & \left[\frac{1}{\hNi}\sum_{k=1}^n|N_i\1_{\{i\}}(\hat I_k)-\hNi\1_{\{i\}}(I_k)|\right]\nonumber \\
& \hspace{3cm}\leq \te|N_i-\hNi|+\te\sum_{k=1}^n|\1_{\{i\}}(\hat I_k)-\1_{\{i\}}(I_k)| \nonumber\\
& \hspace{3cm} \leq 2\te\sum_{k=1}^n|\1_{\{i\}}(\hat I_k)-\1_{\{i\}}(I_k)|.\label{eq:proof-2}
\end{align}
Combining~\eqref{eq:proof-decomp}, \eqref{eq:proof-1} and~\eqref{eq:proof-2} leads to
\begin{align}
\e\bnorme1{(\hfi-\bfi)\1_{N_i>0}} & \leq
2\norme1{K}\sum_{k=1}^n\e\left[\frac{\1_{N_i>0}}{N_i}|\1_{\{i\}}(\hat
  I_k)-\1_{\{i\}}(I_k)|\right] \nonumber \\
& \leq \frac{2\norme1{K}}{n\ai}\sum_{k=1}^n\e\left[\frac{n\ai\1_{N_i>0}}{N_i}|\1_{\{i\}}(\hat I_k)-\1_{\{i\}}(I_k)|\right].\label{eq:proof-3}
\end{align}
The expectation on the right-hand side of this inequality can be bounded in the following way
\begin{align}
\e & \left[\frac{n\ai\1_{N_i>0}}{N_i}|\1_{\{i\}}(\hat
  I_k)-\1_{\{i\}}(I_k)|\right] \leq
\e\left[\frac{n\ai\1_{N_i>0}}{N_i}|\1_{\{i\}}(\hat
  I_k)-\1_{\{i\}}(I_k)|\1_{\frac{n\ai}{N_i}\leq 2}\right] \nonumber\\
&
\hspace{6cm}+\e\left[\frac{n\ai\1_{N_i>0}}{N_i}|\1_{\{i\}}(\hat
  I_k)-\1_{\{i\}}(I_k)|\1_{\frac{n\ai}{N_i}> 2}\right].
\label{eq:proof-4}
\end{align}
For the first term of this bound, we have
\begin{equation}\label{eq:proof-5}
\e\left[\frac{n\ai\1_{N_i>0}}{N_i}|\1_{\{i\}}(\hat I_k)-\1_{\{i\}}(I_k)|\1_{\frac{n\ai}{N_i}\leq 2}\right]\leq
2\varphi_n,
\end{equation}
while for the second term, we obtain from Hölder inequality that
\begin{align}
\e & \left[\frac{n\ai\1_{N_i>0}}{N_i}|\1_{\{i\}}(\hat
  I_k)-\1_{\{i\}}(I_k)|\1_{\frac{n\ai}{N_i}> 2}\right] \nonumber\\
& \leq \sqrt{\e\left[\frac{n\ai\1_{N_i>0}}{N_i}|\1_{\{i\}}(\hat
    I_k)-\1_{\{i\}}(I_k)|\1_{\frac{n\ai}{N_i}>
      2}\right]^2\p\left(\frac{n\ai}{N_i}> 2\right)}\nonumber\\
& \leq
\sqrt{\e\left(\frac{(n\ai)^2}{N_i^2}\1_{N_i>0}\right)
  \p\left(N_i-n\ai<-\frac{n\ai}{2}\right)}.\label{eq:proof-6}
\end{align}
Now, it can be easily seen that
\begin{align}
  \e\left(\frac{(n\ai)^2}{N_i^2}\1_{N_i>0}\right)
  &\leq
  6\e\left(\frac{(n\ai)^2}{(N_i+1)(N_i+2)}\right)\leq 6,\label{eq:proof-7}
\end{align}
where the last inequality follows from~\cite{MR2493852}.
Using  Hoeffding's inequality (see~\cite{hoef63}) we obtain for the second term in~\eqref{eq:proof-4}
\begin{equation}
  \e\left[\frac{n\ai\1_{N_i>0}}{N_i}|\1_{\{i\}}(\hat I_k)-\1_{\{i\}}(I_k)|\1_{\frac{n\ai}{N_i}> 2}\right]\leq \sqrt{6}\exp\left(-\frac{n\ai^2}{4}\right).\label{eq:proof-8}
\end{equation}
From \eqref{eq:proof-3}~--~\eqref{eq:proof-8}, we deduce that
$$\e\bnorme1{(\hat{f}_{i}-\bar{f}_{i})\1_{N_i>0}} \leq\frac{4\norme1{K}}{\ai}\varphi_n+\frac{2\sqrt{6}\norme1K}{\ai}\exp\left(-\frac{n\ai^2}{4}\right).$$
Putting all of the pieces together, we obtain
\begin{align*}
\e\bnorme1{\hfi-\bfi} \leq
\frac{4\norme1{K}}{\ai}\varphi_n+\frac{2\sqrt{6}\norme1K}{\ai} & \exp\left(-\frac{\ai^2}{4}\cdot
n\right)\\
& +\pnorme1K\exp\left(-n\log (1-\ai)\right),
\end{align*}
which concludes the first part of the proof.

Inequality~\eqref{eq:bornealphai} is proved as follows
\begin{align*}
\e|\hai-\ai| & \leq \e\left|\frac{\hNi}{n}-\frac{N_i}{n}\right|+\e\left|\frac{N_i}{n}-\ai\right| \\
& \leq \frac{1}{n}\sum_{k=1}^{n}\e|\1_{\{i\}}(\hat I_k)-\1_{\{i\}}(I_k)|+\frac{1}{n}\sqrt{\var(N_i)} \\
& \leq \varphi_n+\sqrt{\frac{\ai(1-\ai)}{n}}.
\end{align*}

\subsection{Proof of Proposition~\ref{prop:toy_ex}}
Let $k$ be an arbitrary integer in $\interval{n}$. We have to bound $\p(\hat I_k\neq i\given I_k=i)$ for $i=1,2$. To do so, we first consider the case $i=2$:
\begin{align*}
  \p(\hat I_k\neq 2\given I_k=2)
  &= \p(\hat I_k\neq 2, 1-\lambda_n < X_k < 1 \given I_k=2)\\
  &\hspace{2.5cm}+ \p(\hat
  I_k\neq 2, X_k \geq  1 \given I_k=2)\\
  &= \p(1-\lambda_n < X_k < 1 \given I_k=2)
\end{align*}
because, by definition, $\hat I_k \neq 2 \iff X_k < 1$. Thus
\begin{equation}\label{eq:4}
  \p(\hat I_k\neq 2\given I_k=2) = \int_{1-\lambda_n}^1 g_{2,n}(x) dx
  = \lambda_n.
\end{equation}
Next, if $i=1$ it is easy to see that
$\p(\hat I_k\neq 1 \given I_k=1) = \p(X_k\geq 1-\hat\lambda_n \given I_k=1)$.
Let us consider
\begin{equation*}
  \mu_n=\lambda_n+\frac2{\alpha_2}\cdot\frac{\log n}n
  \quad\text{and}\quad
  A=\left\{1-\hat\lambda_n \geq 1-\mu_n\right\}.
\end{equation*}
Using these notations we obtain
\begin{align*}
  \{X_k\geq 1-\lambda_n\} &= \big(\{X_k\geq 1-\hat\lambda_n\}\cap A\big)
  \cup \big\{X_k\geq 1-\hat\lambda_n\}\cap\bar{A}\big)\\
  &\subseteq \{X_k\geq 1-\mu_n\} \cup \left\{\hat\lambda_n\geq\mu_n\right\}.
\end{align*}
This leads to the following inequality
\begin{align}
  \p(\hat I_k\neq 1 \given I_k=1)
  &\leq \mu_n + \p\left(X_{(n)}\leq
    2-\mu_n\given I_k=1 \right)\label{eq:1}.
\end{align}
Since $X_\ell$ and $I_k$ are independent for $k\neq\ell$, we obtain the following bound for the last probability
\begin{align*}
  &\p\left(X_{(n)}\leq
    2-\mu_n\given I_k=1 \right)\\
  &\hspace{2.5cm}= \p\left(\forall \ell, X_{\ell}\leq
    2-\mu_n\given I_k=1 \right)\\
  &\hspace{2.5cm}= \left(\prod_{\ell\neq k}\p\left( X_{\ell}\leq
    2-\mu_n\right)\right)\p\left( X_{k}\leq
    2-\mu_n\given I_k=1 \right).
\end{align*}

The independence of the $X_\ell$'s and simple calculations lead to
\begin{align}
  \p\left(X_{(n)}\leq
    2-\mu_n\given I_k=1 \right) &=
  \big(\p(X_1\leq 2-\mu_n) \big)^{n-1}\nonumber\\
  &= (1-2n^{-1}(\log n))^{n-1}\nonumber\\
  &\leq n^{-1},\label{eq:2}
\end{align}
where the last inequality follows, for $n\geq2$, from the fact that $1-u\leq e^{-u}$ for all $u\geq0$. Taking together equations~\eqref{eq:1} and~\eqref{eq:2}, we finally obtain
\begin{equation}\label{eq:3}
   \p(\hat I_k\neq 1 \given I_k=1) \leq \lambda_n + n^{-1} +
   \frac2{\alpha_2}\cdot\frac{\log n}{n}.
\end{equation}
Proposition follows from equations~\eqref{eq:4} and~\eqref{eq:3}.

\subsection{Proof of Theorem~\ref{theo:phi_suppdisj}}
\label{proof:phi_suppdisj}
Since $\delta_n>2r_n$ we have for all $(i,j)\in\interval{M}^2$ with
$i\neq j$:
\begin{equation}
\label{eq:incsupport}
\left(\bigcup_{k: X_k\in\Sin} B(X_k, r_n)\right)\cap\left(\bigcup_{k:X_k\in S_{j,n}} B(X_k, r_n)\right) \subseteq (\Sin+r_n)\cap (S_{j,n}+r_n)=\emptyset,
\end{equation}
where, for $S\subset\R^d$ and $r>0$, we recall that
$$S+r=\{x\in\R^d:\exists y\in S\textrm{ such that }\|x-y\|_2\leq r\}.$$
Inclusion~\eqref{eq:incsupport} implies $\hat M_{r_n}\geq M$. Moreover, observe that if
\begin{equation}
\label{eq:cond_rn}
r_n\in\mathcal R_M=\{r>0:\hat M_r\leq M\}
\end{equation}
then $\hat M_{r_n}=M$ and the affinity matrices $A^{r_n}$ and $A^{\hat r_n}$ defined in~\eqref{eq:defmatA} induce the same clusters $\cluster1(r_n),\hdots,\cluster{M}(r_n)$.
Furthermore, if \eqref{eq:cond_rn} is verified, it is easily seen that $\forall i\in\interval{M},\exists j\in\interval{M}$ such that
$$\{X_k:X_k\in\Sin+r_n\}\subseteq \cluster{j}(r_n).$$
For simplicity, when \eqref{eq:cond_rn} is satisfied, we index clusters $\cluster1(r_n),\hdots,\cluster{M}(r_n)$ such that
$$\{X_k:X_k\in\Sin+r_n\}\subseteq \cluster{i}(r_n),\quad i\in\interval{M}.$$
We deduce that
\begin{align}
\p(\hat I_k\neq I_k) & \leq \p(\{\hat I_k\neq I_k\}\cap\{r_n\in\mathcal R_M\})+\p(r_n\notin\mathcal R_M)\nonumber \\
& \leq \p(\{\hat I_k\neq I_k\}\cap\{r_n\in\mathcal R_M\}\cap\{X_k\in (S_n+r_n)\}) \nonumber\\
& \hspace{3.5cm}+\p(X_k\notin (S_n+r_n))+\p(r_n\notin\mathcal R_M)\nonumber \\
& \leq \sum_{i=1}^M\p(X_k\notin(\Sin+r_n)|I_k=i)\p(I_k=i)+\psi_n+\p(r_n\notin\mathcal R_M) \nonumber\\
& \leq 2\psi_n+\p(r_n\notin\mathcal R_M)\label{proof:thm2-1}
\end{align}
since $\p(X_k\notin (S_n+r_n))\leq\psi_n$. To complete the proof, we have to find an upper bound for the probability of the event $\left\{r_n\notin\mathcal{R}_M\right\}$. Observe that
\begin{align}
\label{eq:prbarn}
\p(r_n\notin\mathcal R_M)\leq & \p\left(S_n\not\subseteq\bigcup_{k\in\kappa_n}B(X_k,r_n)\right) \nonumber\\
& \hspace{2.5cm}+\p\left(\{r_n\notin\mathcal R_M\}\cap\left\{S_n\subseteq\bigcup_{k\in\kappa_n}B(X_k,r_n)\right\}\right)
\end{align}
where $\kappa_n=\{k\in\interval{M}:X_k\in S_n\}$.
For the first term on the right hand side of the above equation, remark that inclusion
$$S_n\subseteq\bigcup_{k\in\kappa_n}B(X_k,r_n)$$
holds when for all  $\ell\in\interval{N}$, the balls $B_\ell$ defined in assumption 2 contain at least one observation among $\{X_k,k\in\kappa_n\}$. Thus
\begin{align*}
\p\left(S_n\not\subseteq\bigcup_{k\in\kappa_n}B(X_k,r_n)\right) & \leq \p\left(\exists \ell\in\interval{N}, \forall k\in\kappa_n,X_k\notin B_\ell \right) \\
& \leq \sum_{\ell=1}^N\p(\forall k\in\kappa_n,X_k\notin B_\ell) \\
& \leq\sum_{\ell=1}^N\p\left(\bigcap_{k=1}^n\left\{\left\{\{X_k\in S_n\}\cap \{X_k\notin B_\ell\}\right\}\cup\{X_k\notin S_n\}\right\}\right) \\
& \leq \sum_{\ell=1}^N(\p(\{X_k\in S_n\}\cap\{X_k\notin B_\ell\})+\p(X_k\notin S_n))^n \\
& \leq \sum_{\ell=1}^N(1-\p(X_k\in(B_\ell\cap S_n))-\p(X_k\notin S_n)+\p(X_k\notin S_n))^n \\
& \leq \sum_{\ell=1}^N(1-\p(X_k\in(B_\ell\cap S_n)))^n.
\end{align*}
According to assumption 2 and inequality~\eqref{eq:majorant-N}, we obtain
\begin{align*}
\p\left( S_n\not\subseteq\bigcup_{k\in\kappa_n}B(X_k,r_n)\right) & \leq \sum_{\ell=1}^N(1-t_nc_2r_n^d)^n \\
&\leq N \left(1-c_2t_nr_n^d\right)^n\\
&\leq (\tau c_1c_2)^{-1} \frac{n}{\log n} \exp(-c_2 nt_n r_n^d)\\
&\leq (\tau c_1c_2)^{-1} \frac{n}{\log n} \exp(-c_2\tau\log n).
\end{align*}
Since $c_2\tau\geq 1+a$ we have
\begin{equation}\label{proof:thm2-2}
\p\left( S_n\not\subseteq\bigcup_{k\in\kappa_n}B(X_k,r_n)\right)\leq (\tau c_1c_2)^{-1}\frac{1}{n^a\log n}.
\end{equation}
For the second term on the right hand side of \eqref{eq:prbarn}, we have
\begin{equation}\label{proof:thm2-3}
\p\left(\{r_n\notin\mathcal R_M\}\cap\left\{\bigcup_{k\in\kappa_n}B(X_k,r_n)\right\}\right)\leq\p(\exists k\in\interval{n}:X_k\notin (S_n+r_n))\leq n\psi_n.
\end{equation}

Taking \eqref{proof:thm2-1}, \eqref{proof:thm2-2} and \eqref{proof:thm2-3} together, result follows.

\bibliography{akr}
\newpage

\end{document}